\documentclass[12pt]{article}
\usepackage{epsfig}
\usepackage{amsfonts}
\usepackage{amssymb}
\usepackage{amsmath}
\usepackage{amsthm}
\usepackage{latexsym}
\usepackage{graphicx}
\usepackage{bm}
\usepackage{tabularx}
\usepackage{tikz}
\usepackage{booktabs}
\usepackage{enumerate}
\usepackage{caption}
\usepackage[titletoc]{appendix}
\usepackage[numbers]{natbib}
\usepackage{float}
\usepackage[tight,footnotesize]{subfigure}
\usepackage[margin=1in]{geometry}
\usepackage{mwe}
\usepackage{tikz}
\usetikzlibrary{arrows}

\catcode`\@=11

\newcommand{\E}[1]{{\mathrm{E}\left[#1\right]}}

\newcommand{\PP}[1]{\mathrm{P}\left( #1 \right)}
\newcommand{\Var}[1]{{\mathrm{Var}\left(#1\right)}}
\newcommand{\Cov}[1]{{\mathrm{Cov}[#1]}}
\newcommand{\T}{{\mathrm{T}}}

\usepackage[pdftitle={Batch Arrivals},
  pdfauthor={DawPender},
  pdffitwindow=true,
  breaklinks=true,
  colorlinks=true,
  urlcolor=blue,
  linkcolor=red,
  citecolor=red,
  anchorcolor=red]{hyperref}

\numberwithin{equation}{section}
\numberwithin{table}{section}
\numberwithin{figure}{section}

\theoremstyle{plain}
\newtheorem{theorem}{Theorem}[section]
\newtheorem{lemma}[theorem]{Lemma}
\newtheorem{corollary}[theorem]{Corollary}
\newtheorem{proposition}[theorem]{Proposition}

\theoremstyle{definition}

\theoremstyle{remark}
\newtheorem*{remark}{Remark}

\begin{document}

\title{On the Distributions of Infinite Server Queues with Batch Arrivals}
\author{
  Andrew Daw \\ School of Operations Research and Information Engineering \\ Cornell University
\\ 257 Rhodes Hall, Ithaca, NY 14853 \\  amd399@cornell.edu  \\
\and
  Jamol Pender \\ School of Operations Research and Information Engineering \\ Cornell University
\\ 228 Rhodes Hall, Ithaca, NY 14853 \\  jjp274@cornell.edu  \\
 }

\maketitle

\abstract{
Queues that feature multiple entities arriving simultaneously are among the oldest models in queueing theory, and are often referred to as ``batch'' (or, in some cases, ``bulk'') arrival queueing systems. In this work we study the affect of batch arrivals on infinite server queues. We assume that the arrival epochs occur according to a Poisson process, with treatment of both stationary and non-stationary arrival rates. We consider both exponentially and generally distributed service durations and we analyze both fixed and random arrival batch sizes. In addition to deriving the transient mean, variance, and moment generating function for time-varying arrival rates, we also find that the steady-state distribution of the queue is equivalent to the sum of scaled Poisson random variables with rates proportional to the order statistics of its service distribution. We do so through viewing the batch arrival system as a collection of correlated sub-queues. Furthermore, we investigate the limiting behavior of the process through a batch scaling of the queue and through fluid and diffusion limits of the arrival rate. In the course of our analysis, we make important connections between our model and the harmonic numbers, generalized Hermite distributions, and truncated polylogarithms.
}

\section{Introduction}\label{Intro}

Queueing systems with batch arrivals have enjoyed a long and rich history of study, at least on the time scale of queueing theory. Researchers have been exploring models of this sort for no less than six decades, based on the April 1958 submission date of \citet{miller1959contribution}. Given this stretch of time, a wide variety of systems and settings have been considered under the banner of batch arrivals.  Much of the earliest work focuses on single server models, including \citet{miller1959contribution, lucantoni1991new, masuyama2002analysis, liu1993autocorrelations} and \citet{foster1964batched}, although infinite server models followed soon after, such as work by \citet{shanbhag1966infinite} and \citet{brown1969some}. Later work has expanded the concept into a variety of related models, such as for priority queues \cite{takagi1991priority} and for handling server vacations \cite{lee1995batch}. Additionally, there is some work that proves heavy traffic limit theorems for queues with batch arrivals. Examples of this include  \citet{chiamsiri1981diffusion, pang2012infinite, pender2013poisson}.  These papers show that one can approximate the queue length process with Brownian motion and Ornstein-Uhlenbeck processes and also show that one can exploit the approximations even in multi-server and non-Markovian settings.

In this paper we consider queues with arrivals occurring at times following a Poisson process, with consideration given to both non-stationary and stationary rates. We analyze both both general and exponential service as conducted by infinitely many servers. Additionally, this work addresses both fixed and random batch sizes. Our analysis starts with the fixed batch size case. We begin by analyzing the transient behavior of the queue with Markovian service and time-varying arrival rates, providing explicit forms for the moment generating function, mean, and variance. Then, we show that if the arrival rate is stationary the resulting steady-state distribution can be written as a sum of independent, non-identical, scaled Poisson random variables. This leads us to uncover connections to the harmonic numbers and generalizations of the Hermite distribution. By viewing the batch arrival queue as a collection of infinite server sub-queues that receive solitary arrivals simultaneously, we are able to extend this Poisson sum construction to general service distributions. This perspective also provides an avenue for us to extend to random batch sizes. We also give fluid and diffusion scalings of the queue in the case of random batch sizes, as well as extending many of the results we found for fixed batch sizes.

One can note that the batch arrival queue may not always be given the name ``batch,'' as many authors choose to use the term ``bulk'' instead. Predominantly, this reflects two leading strands of applications, where ``bulk'' often gives a connotation of transportation settings whereas ``batch'' frequently implies applications in communications. Just as practical by any other name, this family of models has also been studied in a wide variety of applications beyond these two. Perhaps one most distinct from other types of queueing models is particle splitting in DNA caused by radiation, as discussed in \citet{sachs1992dna}. In this application, primary particles arrive at a cell nucleus and cause DNA double-strand breaks. These double-strand breaks occur in near simultaneity and are thus modeled as arriving in batches of random size, as it is possible that any number double-strand breaks will be induced. After they are induced, the double-strand breaks are then processed by cellular enzymes, corresponding to service in the queueing model. Another interesting and modern application of these models is in cloud-based data processing. In this case, the batches arriving to the system are collections of jobs submitted simultaneously. These jobs are then served by each being processed individually and returned. For more discussion, detailed models, and specific analysis for this setting, see works such as \citet{lu2011join, pender2016law, xie2017pandas, yekkehkhany2018gb} and references therein.
%
%





%

 \subsection{Main Contributions of Paper}

Our contributions in this work can be summarized as follows:

\begin{enumerate}[i)]
\item We show that an infinite server queue with batch arrivals at Poisson process epochs is equivalent in steady state distribution to a sum of scaled independent Poisson random variables, including for generally distributed service and randomly distributed batch sizes. For exponential service, this reveals a connection to the harmonic numbers and generalized Hermite distributions.
\item We derive a limit of the process in which the batch size grows infinitely large and the number of entities in the system is scaled inverse proportionally, yielding a novel distribution characterized by the exponential integral functions. For distributions that meet a divisibility condition, we find that this also holds for random batch sizes.
\item In the case of time-varying arrival rates we give a transient moment generating function for fixed batch sizes as well as means and variance for both fixed and randomly sized batches.
\item We give fluid and diffusion limits of the queue for stationary arrival rates for batches of random size.
\end{enumerate}


\subsection{Organization of Paper}

The body of the remainder of this paper is organized in two main sections: Sections~\ref{deterSec} and~\ref{randomSec}. In Section~\ref{deterSec} we consider systems in which the size of the batches is fixed. Similarly, we devote Section~\ref{randomSec} to the case of randomly distributed batch sizes. At the beginning of each section we give a detailed overview of the contents within and provide context for the analysis in term of this project's scope. After these sections we conclude in Section~\ref{concSec}.


\section{Batches of Deterministic Size}\label{deterSec}

In this section we will consider infinite server queues with arrivals occurring in batches of a fixed size. We will assume that the arrival epochs occur according to a Poisson process, including both stationary and non-stationary models. We also will investigate both exponentially and generally distributed service.

This section starts with studying the case of Markovian arrivals and service in transient state in Subsection~\ref{transientSub}. For a time-varying arrival rate, we give the mean, variance, and moment generating function. We then use this in Subsection~\ref{markovSS} to find the steady-state distribution of the queue. Upon observing that this can be represented as a sum of scaled Poisson random variables, we establish connections to generalized Hermite distributions and to the harmonic numbers. Taking motivation from this, we derive the distribution of the limit of the scaled system as the batch size grows infinitely large. Finally, in Subsection~\ref{orderStatSubsec}, we examine the batch queue as a collection of infinite server sub-queues that simultaneously receive solitary arrivals. In doing so we extend our understanding of the steady-state distribution to the case of general service.

\subsection{Transient Analysis of the Markovian Setting}\label{transientSub}

We begin our analysis with the case of non-stationary Poisson arrival epochs and Markovian service. In Kendall notation, this is the $M_t^n/M/\infty$ queue. We let $Q_t$ represent the number of entities present in the queueing system at time $t \geq 0$, which we often refer to as the ``number in system.'' We will use this notation throughout the remainder of this work, where the precise setting of the queue will be implied by context. In this fully Markovian setting, we can use Dynkin's infinitesimal generator theorem to support our analysis. Specifically, we can note that for a sufficiently regular function $f:\mathbb{N} \to \mathbb{R}$, we have
\begin{align}\label{dynkin}
\frac{\mathrm{d}}{\mathrm{d}t}\E{f(Q_t)}
&
=
\E{\lambda(t) \left(f(Q_t + n) - f(Q_t)\right) + \mu Q_t \left(f(Q_t - 1) - f(Q_t)\right)}
,
\end{align}
for a batch arrival queue with arrival intensity $\lambda(t) > 0$. We will see in this subsection that this infinitesimal generator approach gives us a potent toolkit for exploring this model. Moreover, the insights we find in Markovian settings now and in Subsection~\ref{markovSS} will provide intuition that will guide our investigation of this system when the Markov property does not hold. To begin, we now derive the moment generating function of the number in system. We do so for a system with a non-stationary arrival rate given by a Fourier series, allowing these results to hold for all periodic arrival patterns.

\begin{proposition}\label{batchMGF}
For $\theta \in \mathbb{R}$, let $\mathcal{M}(\theta, t) = \E{e^{\theta Q_t}}$ be the moment generating function of the number in system of an infinite server queue with periodic arrival rate $\lambda + \sum_{k=1}^\infty a_k \cos(k t) + b_k \sin(k t) > 0$, arrival batch size $n \in \mathbb{Z}^+$, and exponential service rate $\mu > 0$. Then, $\mathcal{M}(\theta, t)$ is given by
\begin{align}
\mathcal{M}(\theta, t)
&
=
\left(e^{- \mu t}(e^{\theta} - 1) + 1\right)^{Q_0}
e^{
\sum_{j=1}^n {n \choose j} (e^\theta - 1)^j
\left(
 \frac{\lambda}{j\mu}
 \left(
 1 - e^{-j \mu t}
 \right)
 +
 \sum_{k=1}^\infty
 \frac{(a_k j \mu - b_k k)  }{k^2+j^2\mu^2}
  \left(
 \cos(kt) - e^{-j\mu t}
 \right)
 \right)
 }
 \nonumber
 \\
 &
 \quad
 \cdot
 e^{
\sum_{j=1}^n {n \choose j} (e^\theta - 1)^j
 \sum_{k=1}^\infty
 \frac{(a_k k  +  b_k j \mu )\sin(k t)}{k^2+j^2\mu^2}
}
\end{align}
for all time $t \geq 0$, where $Q_0$ is the initial number in system.
\begin{proof}
From Equation~\ref{dynkin}, the MGF is given by the solution to the partial differential equation
\begin{align*}
\frac{\partial}{\partial t} \mathcal{M}(\theta,t)
&=
\left(\lambda + \sum_{k=1}^\infty a_k \cos(k t) + b_k \sin(k t)\right) \left(e^{n \theta } - 1\right)
\mathcal{M}(\theta,t)
+
\mu \left(e^{-\theta} - 1\right)
\frac{\partial}{\partial \theta}
\mathcal{M}(\theta,t)
\end{align*}
with initial solution $\mathcal{M}(\theta , 0) = e^{\theta Q_0}$. Because $\frac{\mathrm{d}\log(f(x))}{\mathrm{d}x} = \frac{1}{f(x)}\frac{\mathrm{d}f(x)}{\mathrm{d}x}$, we can observe that the partial differential equation for the cumulant generating function $G(\theta, t) = \log \left( \E{e^{\theta Q_t}}\right)$ is
$$
\mu(1 - e^{-\theta})\frac{\partial G(\theta, t)}{\partial \theta}  + \frac{\partial G(\theta, t)}{\partial t} = \left(\lambda + \sum_{k=1}^\infty a_k \cos(k t) + b_k \sin(k t)\right)(e^{n\theta} - 1)
,
$$
with the initial condition
$
G(\theta, 0) = \log\left(\E{e^{\theta Q_0}}\right)  = \theta Q_0
$.
We will now solve this PDE by the method of characteristics. We begin by establishing the characteristic ODE's and corresponding initial solutions as follows:
\begin{align*}
\frac{\mathrm{d}\theta}{\mathrm{d}s}(r,s) &= \mu(1 - e^{-\theta}), &\theta(r,0) = r,\\
\frac{\mathrm{d}t}{\mathrm{d}s}(r,s) &= 1, & t(r,0) = 0,\\
\frac{\mathrm{d}g}{\mathrm{d}s}(r,s) &= \left(\lambda + \sum_{k=1}^\infty a_k \cos(k t) + b_k \sin(k t)\right)(e^{n\theta} - 1), & g(r,0) = rQ_0.
\end{align*}
The first two of these initial value problems yield the following solutions.
\begin{align*}
\theta(r,s) &= \log(e^{c_1(r) + \mu s} + 1) &&\to\quad \theta(r,s) = \log\left((e^r - 1) e^{\mu s} + 1\right)\\
t(r,s) &= s + c_2(r) &&\to\quad t(r,s) = s
\end{align*}
Therefore we can simplify the remaining characteristic ODE to
\begin{align*}
\frac{\mathrm{d}g}{\mathrm{d}s}(r,s)
&
=
\left(\lambda + \sum_{k=1}^\infty a_k \cos(k s) + b_k \sin(k s)\right)\left(\left((e^r - 1) e^{\mu s}+1\right)^n -1 \right)
\\
&
=
\left(\lambda + \sum_{k=1}^\infty a_k \cos(k s) + b_k \sin(k s)\right) \sum_{j=1}^n {n \choose j} (e^r - 1)^j e^{j\mu s},
\end{align*}
and this produces the general solution of
\begin{align*}
g(r,s)
&
=
c_3(r)
+
\sum_{j=1}^n {n \choose j} (e^r - 1)^j
 \left(
 \frac{\lambda}{j\mu}
 +
 \sum_{k=1}^\infty
 \frac{(a_k j \mu - b_k k)  \cos(k s)}{k^2+j^2\mu^2}
 +
 \frac{(a_k k  +  b_k j \mu )\sin(k s)}{k^2+j^2\mu^2}
 \right)
 e^{j\mu s}
 .
 \intertext{This now equates to}
  g(r,s)
&
=
r Q_0
+
\sum_{j=1}^n {n \choose j} (e^r - 1)^j
\Bigg(
 \frac{\lambda}{j\mu}
 \left(
 e^{j \mu s} - 1
 \right)
 +
 \sum_{k=1}^\infty
 \frac{(a_k j \mu - b_k k)  }{k^2+j^2\mu^2}
 \left(
 \cos(ks) e^{j \mu s} - 1
 \right)
 \\
 &
 \quad
 +
 \sum_{k=1}^\infty
 \frac{(a_k k  +  b_k j \mu )\sin(k s)}{k^2+j^2\mu^2}e^{j\mu s}
 \Bigg)
\end{align*}
as the solution to the initial value problem. We now find the solution to the original PDE by solving for each characteristic variable in terms of $t$ and $\theta$ and then substituting these expression into $g(r,s)$. That is, for $s = t$ and $r = \log\left(e^{- \mu t}(e^{\theta} - 1) + 1\right)$, we have that
\begin{align*}
G(\theta, t)
&=
g\left(\log\left(e^{- \mu t}(e^{\theta} - 1) + 1\right), t\right)
\\&=
\log\left(e^{- \mu t}(e^{\theta} - 1) + 1\right)Q_0
+
\sum_{j=1}^n {n \choose j} (e^\theta - 1)^j
\Bigg(
 \frac{\lambda}{j\mu}
 \left(
 1 - e^{-j \mu t}
 \right)
 +
 \sum_{k=1}^\infty
 \frac{(a_k j \mu - b_k k)  }{k^2+j^2\mu^2}
 \\
 &
 \quad
 \cdot
  \left(
 \cos(kt) - e^{-j\mu t}
 \right)
 +
 \sum_{k=1}^\infty
 \frac{(a_k k  +  b_k j \mu )}{k^2+j^2\mu^2} \sin(k t)
 \Bigg)
.
\end{align*}
To conclude the proof, we note that $\mathcal{M}(\theta, t) = e^{G(\theta, t)}$.
\end{proof}
\end{proposition}

We now extend this analysis through two following corollaries. First, for systems with a stationary arrival rate, say $\lambda > 0$, we can further specify the moment generating function explicitly in Corollary~\ref{statCor}. This will be of use when we explore the distribution of the queue in steady-state, which we begin in Subsection~\ref{markovSS}. As with Proposition~\ref{batchMGF}, the uniqueness of moment generating functions will aid us in later exploration of the distributions within this model and within generalizations of it.

\begin{corollary}\label{statCor}
For $\theta \in \mathbb{R}$, let $\mathcal{M}(\theta, t) = \E{e^{\theta Q_t}}$ be the moment generating function of the number in system of an infinite server queue with stationary arrival rate $\lambda > 0$, arrival batch size $n \in \mathbb{Z}^+$, and exponential service rate $\mu > 0$. Then, $\mathcal{M}(\theta, t)$ is given by
\begin{align}
\mathcal{M}(\theta, t)
&
=
\left(e^{- \mu t}(e^{\theta} - 1) + 1\right)^{Q_0}
e^{
\lambda \sum_{j=1}^n {n \choose j} \frac{(e^\theta - 1)^j }{j\mu}
 \left(
 1 - e^{-j \mu t}
 \right)
 }
\end{align}
for all time $t \geq 0$, where $Q_0$ is the initial number in system.
\end{corollary}

For the second direct result of Proposition~\ref{batchMGF}, we can also give explicit expressions for the transient mean and variance of the queue. We derive these equations from the first and second derivatives, respectively, of the cumulant generating function $\log(\E{e^{Q_t}})$.

\begin{corollary}\label{meanvarcor}
Let $Q_t$ be an infinite server queue with periodic arrival rate $\lambda + \sum_{k=1}^\infty a_k \cos(k t) + b_k \sin(k t) > 0$, arrival batch size $n \in \mathbb{Z}^+$, and exponential service rate $\mu > 0$. Then, the mean and variance of the queue are given by
\begin{align}
\E{Q_t}
&
=
Q_0 e^{-\mu t}
+
 \frac{n\lambda }{\mu}
 \left(
 1 - e^{- \mu t}
 \right)
 +
 \sum_{k=1}^\infty
 \frac{n(a_k  \mu - b_k k)  }{k^2+\mu^2}
  \left(
 \cos(kt) - e^{-\mu t}
 \right)
 \nonumber
 \\
 &
 \quad
 +
 \sum_{k=1}^\infty
 \frac{n(a_k k  +  b_k  \mu )}{k^2+\mu^2} \sin(k t)
\\
\Var{Q_t}
&
=
Q_0\left(e^{-\mu t} - e^{-2\mu t}\right)
+
 \frac{n\lambda }{\mu}
 \left(
 1 - e^{- \mu t}
 \right)
 +
 \sum_{k=1}^\infty
 \frac{n(a_k  \mu - b_k k)  }{k^2+\mu^2}
  \left(
 \cos(kt) - e^{-\mu t}
 \right)
 \nonumber
 \\
 &
 \quad
 +
 \sum_{k=1}^\infty
 \frac{n(a_k k  +  b_k  \mu )}{k^2+\mu^2} \sin(k t)
 +
 \frac{n(n-1) \lambda }{2\mu}
 \left(
 1 - e^{-2 \mu t}
 \right)
 +
 \sum_{k=1}^\infty
 \frac{n(n-1)(2 a_k \mu - b_k k)  }{k^2 + 4\mu^2}
 \nonumber
 \\
 &
 \quad
 \cdot
  \left(
 \cos(kt) - e^{-2\mu t}
 \right)
 +
 \sum_{k=1}^\infty
 \frac{n(n-1)(a_k k  +  2 b_k \mu )}{k^2+ 4\mu^2} \sin(k t)
\end{align}
for all time $t \geq 0$, where $Q_0$ is the initial number in system.
\end{corollary}


In the remainder of this work we will explore various modifications of this model, including general service and randomized batch sizes. The results of this subsection will serve as cornerstone throughout much of this upcoming analysis, both supporting the underlying derivation techniques and providing the intuition for new perspectives.

\subsection{The Markovian System with Stationary Arrival Rates}\label{markovSS}

Our first departure from our initial model will be modest: instead of studying the fully Markovian, non-stationary, fixed batch size system in transient time we will now move to addressing the stationary case, with much of our analysis focused on the system in steady-state. This simplified setting will allow us to extract greater intuition from our prior findings, which in turn will support generalization of the service distribution and randomization of the batch sizes. To begin, we find a representation of the steady-state distribution of the queue length in terms of a sum of independent, scaled Poisson random variables.

\begin{proposition}\label{ssDist}
In steady-state the distribution of the number in system of an infinite server queue with stationary arrival rate $\lambda > 0$, arrival batch size $n \in \mathbb{Z}^+$, and exponential service rate $\mu > 0$ is
\begin{align}
Q_\infty(n)
\stackrel{D}{=}
\sum_{j=1}^n
j Y_j
\end{align}
where $Y_j \sim \mathrm{Pois}\left(\frac{\lambda}{j\mu}\right)$ are independent.
\begin{proof}
From Proposition~\ref{batchMGF}, we have that the steady-state moment generating function of the queue is given by
\begin{align*}
\lim_{t\to \infty}\mathcal{M}(\theta, t)
&
=
e^{
\lambda \sum_{k=1}^n {n \choose k} \frac{\left(e^{\theta} - 1\right)^k}{k\mu}
}
.
\end{align*}
To satisfy our stated Poisson form, we are now left to show that $\sum_{k=1}^n {n \choose k}\frac{(e^\theta - 1)^k}{k} = \sum_{k=1}^n \frac{e^{k \theta} - 1}{k}$ for all $n \in \mathbb{Z}^+$. We proceed by induction. In the base case of $n=1$ we have $e^{\theta} - 1 = e^\theta - 1$ and so we are left to show the inductive step. We now assume $\sum_{k=1}^n {n \choose k}\frac{(e^\theta - 1)^k}{k} = \sum_{k=1}^n \frac{e^{k \theta} - 1}{k}$ holds at $n$. Then, by the Pascal triangle identity ${n  \choose k} = {n+1 \choose k} - {n \choose k -1}$ and our inductive hypothesis we can observe
\begin{align*}
\sum_{k=1}^n \frac{e^{k \theta} - 1}{k}
&=
\sum_{k=1}^n {n \choose k}\frac{(e^\theta - 1)^k}{k}
=
\sum_{k=1}^n \left({n+1 \choose k} - {n \choose k -1}\right)\frac{(e^\theta - 1)^k}{k}.
\end{align*}
Now, by applying the identity ${n \choose k -1} = \frac{k}{n+1}{n+1 \choose k}$ and distributing the summation we can further note that
\begin{align*}
\sum_{k=1}^n \left({n+1 \choose k} - {n \choose k -1}\right)\frac{(e^\theta - 1)^k}{k}
&
=
\sum_{k=1}^n \left({n+1 \choose k} - \frac{k}{n+1}{n+1 \choose k}\right)\frac{(e^\theta - 1)^k}{k}
\\
&
=
\sum_{k=1}^n {n+1 \choose k}\frac{(e^\theta - 1)^k}{k} - \frac{\sum_{k=1}^n{n+1 \choose k}(e^\theta - 1)^k}{n+1}.
\end{align*}
Now, we can use the binomial theorem to see that
$$
\sum_{k=1}^n{n+1 \choose k}(e^\theta - 1)^k
=
(e^\theta - 1 + 1)^{n+1} - 1 - (e^{\theta}-1)^{n+1}
=
e^{(n+1)\theta} - 1 - (e^{\theta}-1)^{n+1}
,
$$
and so we can now simplify and find
\begin{align*}
\sum_{k=1}^n {n+1 \choose k}\frac{(e^\theta - 1)^k}{k}
-
\frac{\sum_{k=1}^n{n+1 \choose k}(e^\theta - 1)^k}{n+1}
&
=
\sum_{k=1}^n {n+1 \choose k}\frac{(e^\theta - 1)^k}{k}
+
\frac{(e^\theta - 1)^{n+1}}{n+1}
\\
&
\quad
-
\frac{e^{(n+1)\theta} - 1}{n+1}
.
\end{align*}
Hence, in conjunction with our initial equation, we have that
\begin{align*}
\sum_{k=1}^n \frac{e^{k \theta} - 1}{k}
&
=
\sum_{k=1}^n {n+1 \choose k}\frac{(e^\theta - 1)^k}{k}
+
\frac{(e^\theta - 1)^{n+1}}{n+1}
-
\frac{e^{(n+1)\theta} - 1}{n+1}
,
\end{align*}
and by rearranging terms we now complete the inductive approach:
\begin{align*}
\sum_{k=1}^{n+1} \frac{e^{k \theta} - 1}{k}
&
=
\sum_{k=1}^{n+1} {n+1 \choose k}\frac{(e^\theta - 1)^k}{k}
.
\end{align*}
We can now observe that we have a moment generating function that is a product of moment generating functions of scaled Poisson random variables, which yields the stated result.
\end{proof}
\end{proposition}

While we will continue to explore the stationary arrival rate setting throughout this subsection, we note that this Poisson sum representation will be a leading inspiration in the sequel. Specifically, in Subsection~\ref{orderStatSubsec} we will find intuition for this result by viewing the batch arrival queue as a collection of sub-systems.

\begin{remark}
In addition to this Poisson sum representation, we can also express the steady-state MGF in terms of the truncated polylogarithm function and harmonic numbers.  From the MGF of the queue length in steady state for $\theta < 0$, we can observe that
\begin{align*}
\lim_{t\to \infty}\mathcal{M}(\theta, t)
&
=
 e^{ \frac{\lambda}{\mu} \sum_{k=1}^n \frac{e^{k \theta} - 1}{k} }
= e^{ \frac{\lambda}{\mu} ( \mathrm{Li}(e^\theta,n,1) - H_n )}
\end{align*}
where we have  $H_n$ as the $n^\text{th}$ harmonic number, given by $\sum^{n}_{k=1} \frac{1}{k}$, and where the truncated polylogarithm function $\mathrm{Li}(z,n,s)$ is defined as
\begin{equation*}
\mathrm{Li}(z,n,s) = \sum^{n}_{k=1} \frac{z^k}{k^s}.
\end{equation*}
\end{remark}

 This decomposition into Poisson random variables can be quite useful from a computational standpoint.  It allows us to simulate the steady state quite easily since we only need to simulate $n$ Poisson random variables instead of simulating an actual queue, which could be quite expensive. We can now observe that this construction also yields an interesting connection to both the harmonic number and Hermite distributions, as suggested in the remark above. To motivate our following analysis, suppose that $n=2$. Then, steady-state queue length has steady-state moment generating function given by
$$
\mathcal{M}_n(\theta, \infty)
=
e^{ \frac{\lambda}{\mu} \left(e^{ \theta} - 1\right) +  \frac{\lambda}{2\mu} \left(e^{ 2\theta} - 1\right) }.
$$
We can now observe that this MGF corresponds to a Hermite distribution with parameters $\frac{\lambda}{\mu}$ and $\frac{\lambda}{2\mu}$. This implies that the steady-state CDF of the queue at $n=2$ is
\begin{align*}
P(Q_{\infty}(2) \leq k)
& = e^{-\frac{3\lambda}{2\mu} }\sum_{i=0}^{\lfloor k\rfloor} \sum_{j=0}^{\lfloor i/2\rfloor} \frac{\left(\frac{\lambda}{\mu}\right)^{i-2j}\left(\frac{\lambda}{2\mu}\right)^j}{(i-2j)!j!}
 = e^{-\frac{3\lambda}{2\mu} } \sum_{i=0}^{\lfloor k\rfloor} \sum_{j=0}^{\lfloor i/2\rfloor} \frac{\left( \frac{\lambda}{\mu} \right)^{i-j} 2 ^{-j}}{(i-2j)!j!}
.
\end{align*}
Furthermore, the steady-state PMF of the queue length is given by
\begin{align*}
P(Q_{\infty}(2) = i ) & = e^{-\frac{3\lambda}{2\mu} } \sum_{j=0}^{\lfloor i/2\rfloor} \frac{\left( \frac{\lambda}{\mu} \right)^{i-j} 2 ^{-j}}{(i-2j)!j!} .
\end{align*}
This observation prompts us to ponder generalizations for $n \geq 3$. The term ``generalized Hermite distribution'' has taken on slightly varying (yet always interesting) definitions for different authors. For readers interested in the Hermite distribution and popular generalizations of it, we suggest \citet{kemp1965some,gupta1974generalized}, and \citet{milne1993generalized}. In our setting we note that the coefficients of $\frac{\lambda}{\mu}$ in the MGF for batch size $n$ will be $1$, $\frac{1}{2}$, $\frac{1}{3}$, \dots, $\frac{1}{n}$. For this reason, we think of this particular generalization of Hermite distributions to be the \textit{harmonic Hermite distribution}. We can now note that because of this harmonic structure we can instead fully characterize the distribution simply by $n$ and $\frac{\lambda}{\mu}$. In the following proposition we find a useful recursion for the probability mass function of this distribution at all $n \in \mathbb{Z}^+$.

\begin{proposition}\label{recurProp}
Let $Q_{t}(n)$ be an infinite server batch arrivals queue with arrival rate $\lambda > 0$, batch size $n \in \mathbb{Z}^+$, and service rate $\mu > 0$. Then, the steady-state distribution of the queue is given by the recursion
\begin{equation}
\mathbb{P}(Q_{\infty}(n) = j ) = p_j = \sum^{n}_{i=1} i p_{j-i} \frac{\lambda}{i j \mu} = \sum^{n}_{i=1} p_{j-i} \frac{\lambda}{ j \mu},
\end{equation}
where $p_0 = e^{-\frac{\lambda}{\mu}H_n}$ for $H_n$ as the $n^\text{th}$ harmonic number and $p_k = 0$ for all $k < 0$. Thus, we say that $Q_\infty(n)$ follows the ``harmonic Hermite distribution'' with parameter $n$.
\begin{proof}
We know from our Poisson representation of the steady state queue length that the steady-state moment generating function is
\begin{equation*}
M(\theta) = \sum^{\infty}_{j=0} \mathbb{P}(Q_{\infty}(n) = j ) \theta^j = \sum^{\infty}_{j=0} p_j  \theta^j = \exp \left( \sum^{n}_{i=1} \frac{\lambda}{i \mu} \left( \theta^i  - 1 \right)  \right) .
\end{equation*}
If we take the logarithm of both sides we see that we have
\begin{equation*}
\log \left( \sum^{\infty}_{j=0} p_j  \theta^j \right) = \sum^{n}_{i=1} \frac{\lambda}{i \mu} \left( \theta^i  - 1 \right) .
\end{equation*}
Now we take the derivative of both sides with respect to the parameter $\theta$ and this yields the following expression
\begin{equation*}
\frac{ \sum^{\infty}_{j=1} j p_j  \theta^{j-1} }{ \sum^{\infty}_{j=0} p_j  \theta^j } = \sum^{n}_{i=1} \frac{\lambda}{ \mu} \theta^{i-1}  .
\end{equation*}
By moving the denominator to the righthand side, we have that
\begin{equation*}
\sum^{\infty}_{j=1} j p_j  \theta^{j-1} =  \left( \sum^{\infty}_{j=0} p_j  \theta^j  \right) \left(  \sum^{n}_{i=1} \frac{\lambda}{ \mu} \theta^{i-1}  \right).
\end{equation*}
Finally, by matching similar powers of $\theta$ on the left and right sides, we complete the proof.
\end{proof}
\end{proposition}

From the above result, we see that for the steady state queue length $Q_{\infty}(n)$ we can derive the specific probabilities,
\begin{eqnarray*}
p_0 &=& e^{-\frac{\lambda}{\mu}H_n}, \\
p_1 &=&  \frac{\lambda}{\mu} p_0 = \frac{\lambda}{\mu} e^{-\frac{\lambda}{\mu}H_n},  \\
p_2 &=&  \frac{\lambda}{2\mu} ( p_0 + p_1 ) = \frac{\lambda}{2\mu} e^{-\frac{\lambda}{\mu}H_n}  + \frac{\lambda^2}{2\mu^2} e^{-\frac{\lambda}{\mu}H_n} .
\end{eqnarray*}
We can repeat this process as needed for any desired probability. From Proposition~\ref{ssDist}, we can observe that the mean number in system grows linearly with the batch size, meaning that the mean of the $n^\text{th}$ harmonic Hermite distribution is
\begin{align}\label{mean1}
\E{Q_\infty(n)} = \sum_{j=1}^n j \E{Y_j} = \frac{n\lambda}{\mu}.
\end{align}
We can observe further that the second moment and variance are quadratic functions of $n$:
$$
\E{Q_\infty(n)^2} = \E{ \left( \sum_{j=1}^n j Y_j \right)^2} = \frac{n (n+1)\lambda}{2\mu} + n^2 \frac{\lambda^2}{\mu^2},
$$
$$
\mathrm{Var}[Q_\infty(n)] = \E{Q_\infty(n)^2}  -  \E{Q_\infty(n)}^2 = \frac{n (n+1)\lambda}{2\mu} .
$$
We note that from Proposition~\ref{ssDist} and the following remark, the moment generating function of this distribution is given by
\begin{align}\label{mgf1}
\lim_{t\to \infty}\mathcal{M}(\theta, t)
&
=
 e^{ \frac{\lambda}{\mu} \sum_{k=1}^n \frac{e^{k \theta} - 1}{k} }
= e^{ \frac{\lambda}{\mu} ( \mathrm{Li}(e^\theta,n,1) - H_n )}.
\end{align}
If one is to consider this system as the batch size grows infinitely large we can see from Equations~\ref{mean1} and~\ref{mgf1} that the number in system will grow proportionally, tending to infinity as $n$ does. This leads us to ponder the limiting object of the scaled number in system $\frac{Q_t(n)}{n}$ as the batch size grows.

We begin by using Equation~\ref{mgf1} with $\theta$ replaced by $\frac{\theta}{n}$ to see that the steady-state moment generating function of this scaled queue length is
\begin{align}\label{scaleSS}
\lim_{t\to \infty}\mathcal{M}(\theta, t)
=
 e^{ \frac{\lambda}{\mu} \sum_{k=1}^n \frac{e^{\frac{k}{n} \theta} - 1}{k} }
 .
\end{align}
Furthermore, by replacing $\theta$ with $\frac{\theta}{n}$ and $Q_0(n)$ with $\frac{Q_0(n)}{n}$ in Proposition~\ref{batchMGF}, we can note that the transient moment generating function for this scaled system with constant arrival rate is given by
$$
\E{e^{\theta \cdot \frac{Q_{t}(n)}{n}}}  \equiv \mathcal{M}_n(\theta, t)
=
\left(e^{- \mu t}(e^{\frac{\theta}{n}} - 1) + 1\right)^{\frac{Q_0}{n}}
e^{
\lambda \sum_{k=1}^n {n \choose k} \frac{\left(e^{\theta/n} - 1\right)^k}{k\mu} \left(1-e^{-k\mu t}\right)
}.
$$
Additionally, we can also observe that the steady-state distribution of the scaled queue can also be interpreted as a sum of Poisson random variables through direction application of Proposition~\ref{ssDist} or by inspection of Equation~\ref{scaleSS}. This representation is
\begin{align}\label{scalePois}
\frac{Q_\infty(n)}{n}
\stackrel{D}{=}
\sum_{j=1}^n \frac{j}{n} Y_j
,
\end{align}
where again $Y_j \sim \text{Pois}\left(\frac{\lambda}{j \mu}\right)$.

We now consider the limit as $n \to \infty$, in which we are both sending the size of batches of arrivals to infinity while also scaling the size of the queue inversely. We can use this construction to move beyond just the mean and variance and instead explicitly state every cumulant of the scaled queue. In Proposition~\ref{cumulants} we give exact expressions of all steady-state cumulants of the scaled queue as functions of the Bernoulli numbers. Further, we find a convenient form of every cumulant of the scaled queue as the batch size grows to infinity.

\begin{proposition}\label{cumulants}
Let $\lambda > 0$ be the arrival rate of batches of size $n \in \mathbb{Z}^+$ to an infinite server queue with exponential service rate $\mu > 0$. Then, the $k^\text{th}$ steady-state cumulant of the scaled queue $\mathcal{C}^k\left[ \frac{Q_\infty(n)}{n} \right]$ is given by
 \begin{eqnarray}
 \mathcal{C}^k\left[ \frac{Q_\infty(n)}{n} \right] &=& \frac{\frac{n^{k}}{k}+\frac{1}{2}n^{k-1}+\sum_{j=2}^{k-1} \frac{\mathrm{B}_{j}}{j!}(k-1)_{j-1}n^{k-j}}{n^k}.
\end{eqnarray}
where $(n)_{i} = \frac{n!}{(n-i)!} $ is the $i^\text{th}$ falling factorial of $n$ and $\mathrm{B}_i$ is the $i^\text{th}$ Bernoulli number, which is defined as
$$
B_i
=
\sum_{k=0}^i
\sum_{j=0}^k
(-1)^j
{k \choose j}
\frac{(j+1)^i}{k+1}
.
$$
Moreover, we have that
$
\lim_{n \to \infty} \mathcal{C}^k\left[ \frac{Q_\infty(n)}{n} \right] = \frac{\lambda}{k\mu}.
$

\begin{proof}
From our prior observation that
$
\frac{Q_\infty(n)}{n}
\stackrel{D}{=}
\sum_{j=1}^n
\frac{j}{n} Y_j
$
where $Y_j \sim \mathrm{Pois}\left(\frac{\lambda}{j\mu}\right)$, we have that
\begin{eqnarray*}
\mathcal{C}^k\left[ \frac{Q_\infty(n)}{n} \right]
=
  \mathcal{C}^k\left[  \sum_{j=1}^n \frac{j}{n} Y_j \right]
=
\sum_{j=1}^n  \mathcal{C}^k\left[ \frac{j}{n} Y_j \right]
=
\sum_{j=1}^n  \frac{j^k}{n^k} \mathcal{C}^k\left[  Y_j \right]
=
\frac{\lambda}{\mu n^k}  \sum_{j=1}^n  j^{k-1},
\end{eqnarray*}
from the independence of these Poisson distributions.
Now, by using Faulhaber's formula as given in \citet{knuth1993johann}, we achieve the stated result.
\end{proof}
\end{proposition}

Just as we built from inherited expressions for the mean and variance to specify every cumulant in Proposition~\ref{cumulants}, we can also find the limit of the transient-state moment generating function for the scaled queue given in Equation~\ref{mgf1}.

\begin{proposition}\label{scaleMGF}
Let $Q_t$ be an infinite server queue with arrival rate $\lambda > 0$, arrival batch size $n \in \mathbb{Z}^+$, and exponential service rate $\mu > 0$. For $\theta \in \mathbb{R}$, let
$$
\mathcal{M}_\infty(\theta, t) = \lim_{n \to \infty} \E{e^{\frac{\theta Q_t(n)}{n}}}.
$$
Then, $\mathcal{M}_\infty(\theta, t)$ is given by
\begin{align}
\mathcal{M}_\infty(\theta, t)
&
=
\begin{cases}
e^{
\frac{\lambda}{\mu}
 \left(
\mathrm{Ei}(\theta)
-
\mathrm{Ei}(\theta e^{-\mu t} )
-
\mu t
 \right)
}
& \text{ if $\theta > 0$,}
\\
e^{
 \frac{\lambda}{\mu}
 \left(
 E_1(-\theta e^{-\mu t})
 -
E_1(-\theta)
-
\mu t
\right)
}
& \text{ if $\theta < 0$,}
\\
1
& \text{ if $\theta = 0$,}
\end{cases}
\end{align}
for all time $t \geq 0$, where the exponential integral functions $\mathrm{Ei}(x)$ and $E_1(x)$ are defined
$$
\mathrm{Ei}(x) = -\int_{-x}^\infty \frac{e^{-s}}{s} \mathrm{d}s,
\quad
E_1(x) = \int_{x}^\infty \frac{e^{-s}}{s} \mathrm{d}s,
$$
and are real-valued for $x > 0$.
\begin{proof}
While conventions may vary by application area, in this work we use the definition of exponential integral function given by
$$
\mathrm{Ei}(x) = -\int_{-x}^\infty \frac{e^{-s}}{s} \mathrm{d}s.
$$
By taking the limit of the MGF of the scaled queue, we have that
\begin{align*}
\frac{\partial}{\partial t} \mathcal{M}_\infty(\theta,t)
&=
\lambda \left(e^{ \theta } - 1\right)
\mathcal{M}_\infty(\theta,t)
-
\mu \theta \frac{\partial}{\partial \theta}
\mathcal{M}_\infty(\theta,t)
\end{align*}
with initial solution $\mathcal{M}_\infty(\theta , 0) = \lim_{n \to \infty} e^{\frac{\theta Q_0}{n}} = 1$. In the same manner as the proof of Theorem~\ref{batchMGF}, we solve the PDE of the cumulant generating function through use of the method of characteristics. We start by establishing the characteristic ODE's:
\begin{align*}
\frac{\mathrm{d}\theta}{\mathrm{d}s}(r,s) &= \mu\theta , &\theta(r,0) = r,\\
\frac{\mathrm{d}t}{\mathrm{d}s}(r,s) &= 1, & t(r,0) = 0,\\
\frac{\mathrm{d}g}{\mathrm{d}s}(r,s) &= \lambda(e^{\theta} - 1), & g(r,0) = 0.
\end{align*}
We now solve the first two initial value problems and find
\begin{align*}
\theta(r,s) &= c_1(r)e^{\mu s} &&\to\quad \theta(r,s) = r e^{\mu s},\\
t(r,s) &= s + c_2(r) &&\to\quad t(r,s) = s.
\end{align*}
This allows us to simplify the third characteristic equation to
$$
\frac{\mathrm{d}g}{\mathrm{d}s}(r,s)
=
\lambda(e^{r e^{\mu s}} - 1).
$$
Because $\theta = r e^{\mu s}$, we can note that $r$ and $\theta$ will match in sign: $r > 0$ if and only if $\theta > 0$. If $\theta > 0$, the general solution to this ODE is
\begin{align*}
g(r,s)
&
=
 c_3(r)
 +
 \frac{\lambda}{\mu}
 \left(
\mathrm{Ei}(r e^{\mu s})
-
\mu s
 \right)
 ,
\end{align*}
whereas if $\theta < 0$, the solution is instead
\begin{align*}
g(r,s)
&
=
 c_3(r)
 -
 \frac{\lambda}{\mu}
 \left(
E_1(-r e^{\mu s})
+
\mu s
 \right)
 .
\end{align*}
This follows from the fact that for $x > 0$ the exponential integral functions are such that $\mathrm{Ei}(x) = -E_1(-x) - i \pi$; that is, the real parts of $E_1(-x)$ and $-\mathrm{Ei}(x)$ are the same. Moreover, for $x > 0$ one can consider $\mathrm{Ei}(x)$ as the real part of $-E_1(-x)$. Additionally, $E_1(x)$ is real for all $x > 0$. Hence, we use each definition of the exponential integral function when appropriate. As an alternative, we could replace each of these functions with $\texttt{real}(-E_1(-x))$ to have a single expression for both positive and negative $x$. For a collection of facts regarding the exponential integral functions, see Pages 228-237 of \citet{abramowitz1965handbook}.

Now, using this we have that the corresponding solutions to the initial value problems will be
\begin{align*}
  g(r,s)
  &
=
\begin{cases}
\frac{\lambda}{\mu}
 \left(
\mathrm{Ei}(r e^{\mu s})
-
\mathrm{Ei}(r )
-
\mu s
 \right)
 &
 \text{ if $r > 0$,}
 \\
 \frac{\lambda}{\mu}
 \left(
 E_1(-r)
 -
E_1(-r e^{\mu s})
-
\mu s
 \right)
 &
 \text{ if $r < 0$}
.
\end{cases}
\end{align*}
Hence, for $s = t$ and $r = \theta e^{-\mu t}$, this yields
\begin{align*}
G(\theta, t)
&=
g\left(\theta e^{-\mu t}, t\right)
=
\begin{cases}
\frac{\lambda}{\mu}
 \left(
\mathrm{Ei}(\theta)
-
\mathrm{Ei}(\theta e^{-\mu t} )
-
\mu t
 \right)
 &
 \text{ if $\theta > 0$,}
 \\
 \frac{\lambda}{\mu}
 \left(
 E_1(-\theta e^{-\mu t})
 -
E_1(-\theta)
-
\mu t
 \right)
 &
 \text{ if $\theta < 0$}
.
\end{cases}
\end{align*}
By $\mathcal{M}_\infty(\theta, t) = e^{G_\infty(\theta, t)}$, we complete the proof.
\end{proof}
\end{proposition}

By consequence, we can also give the moment generating function in steady-state.

\begin{corollary}\label{ssMGF}
The moment generating function of the scaled number in system in steady-state as $n \to \infty$ is given by
\begin{align}
\mathcal{M}_\infty(\theta)
&
=
\begin{cases}
\theta^{-\frac{\lambda}{\mu}}
e^{
\frac{\lambda}{\mu}
 \left(
\mathrm{Ei}(\theta)
-
\gamma
 \right)
}
& \text{ if $\theta > 0$,}
\\
(-\theta)^{-\frac{\lambda}{\mu}}
e^{
-
 \frac{\lambda}{\mu}
 \left(
E_1(-\theta)
+
\gamma
\right)
}
& \text{ if $\theta < 0$,}
\\
1
& \text{ if $\theta = 0$,}
\end{cases}
\end{align}
where $\gamma$ is the Euler-Mascheroni constant.
\begin{proof}
From \citet{abramowitz1965handbook}, for $x > 0$ we can expand the exponential integral functions as
\begin{align}
\mathrm{Ei}(x) &= \gamma + \log(x) + \sum_{k=1}^\infty \frac{x^k}{k k!} ,
\quad
E_1(x) = - \gamma - \log(x) - \sum_{k=1}^\infty \frac{(-x)^k}{k k!} , \label{e1expand}
\end{align}
where $\gamma$ is the Euler-Mascheroni constant. By expanding $\mathrm{Ei}(\theta e^{-\mu t})$ and $E_1(-\theta e^{-\mu t})$ in the respective cases of positive and negative $\theta$ and taking the limit as $t \to \infty$, we achieve the stated result.
\end{proof}
\end{corollary}

\begin{figure}[h]
\begin{center}	
\includegraphics[width=.9\textwidth]{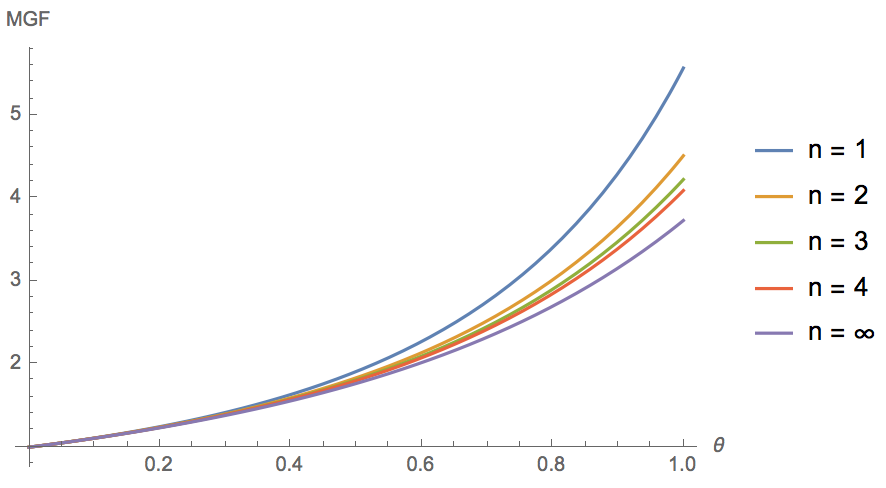}
\caption{Steady-state MGF of the scaled queue for increasing batch size where $\frac{\lambda}{\mu} = 1$.} \label{limMGFscale}
\end{center}
\end{figure}

As a demonstration of the convergence of the steady-state moment generating functions of the batch scaled queues to the expression given in Corollary~\ref{ssMGF}, we plot the first four cases in comparison to the limiting scenario in Figure~\ref{limMGFscale}.

While it can be argued that even in steady-state the form of this moment generating function is unfamiliar, we can still observe interesting characteristics of it. In particular, for $\theta < 0$ we can uncover a connection back to the harmonic numbers. We now discuss this in the following remark.

\begin{remark}
Using Equation~\ref{e1expand}, we can note that for $\theta < 0$ the steady-state moment generating function of limit of the scaled queue can be expressed
$$
M(\theta)
=
(-\theta)^{-\frac{\lambda}{\mu}}
e^{
-
 \frac{\lambda}{\mu}
 \left(
E_1(-\theta)
+
\gamma
\right)
}
=
e^{
-\frac{\lambda}{\mu}\left(
E_1(-\theta)
+
\gamma
+
\log(-\theta)
\right)
}
=
e^{
-\frac{\lambda}{\mu}\left(
-
\sum_{k=1}^\infty \frac{\theta^k}{k k!}
\right)
}
.
$$
From \citet{dattoli2008note}, we have that $-e^{x} \sum_{k=1}^\infty \frac{(-x)^k}{k k!}$ is an exponential generating function for the harmonic numbers. That is,
$$
-e^{x} \sum_{k=1}^\infty \frac{(-x)^k}{k k!}
=
\sum_{n=1}^\infty
\frac{x^n}{n!}H_n
$$
where $H_n$ is the $n^\text{th}$ harmonic number. Thus, for $\theta < 0$ the steady-state moment generating function of this limiting object can be further simplified to
$$
M(\theta)
=
e^{
-\frac{\lambda}{\mu}\left(
-
\sum_{k=1}^\infty \frac{\theta^k}{k k!}
\right)
}
=
e^{
-\frac{\lambda}{\mu}
\sum_{n=1}^\infty
H_n
e^{\theta}
\frac{(-\theta)^n}{n!}
}
=
e^{
-\frac{\lambda}{\mu}
\E{
H_N
}
}
,
$$
where $N \sim \mathrm{Pois}(-\theta)$.
\end{remark}

\begin{figure}[h]
\begin{center}	
\includegraphics[width=.5\textwidth]{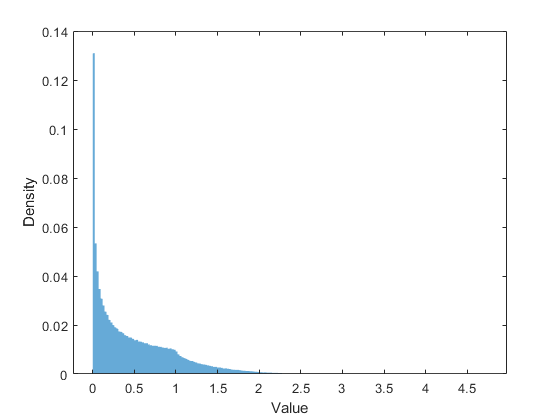}\\
\includegraphics[width=.5\textwidth]{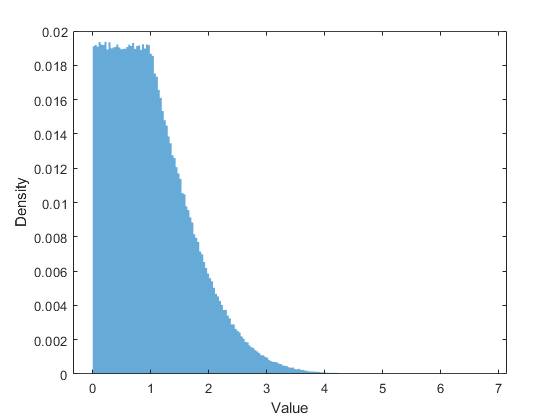}~\includegraphics[width=.5\textwidth]{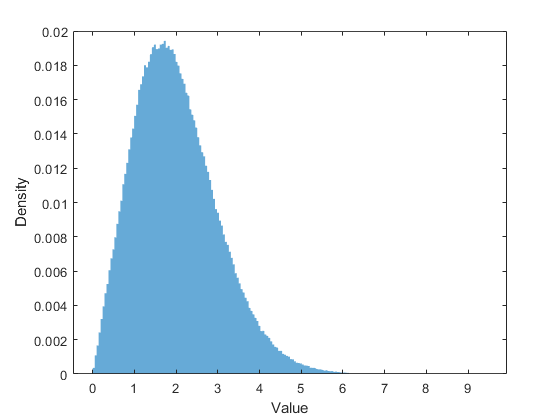}
\caption{Approximate steady-state density of the scaled queue limit for size where $\frac{\lambda}{\mu} = \frac{1}{2}$ (top), $\frac{\lambda}{\mu} = 1$ (left), and $\frac{\lambda}{\mu} = 2$ (right), using 1,000,000 simulation replications and $n=2,000$.} \label{batchLim1}
\end{center}
\end{figure}

In addition to this remark's connection of the moment generating function and the harmonic numbers, we can also gain insight into this limiting object through Monte Carlo methods. Using Equation~\ref{scalePois}, we have a simple and efficient approximate simulation method for this process through summing scaled Poisson random numbers. Furthermore, this approximation of course becomes increasingly precise as $n$ grows. As an example of this, we give the simulated steady-state densities across different relationships of $\lambda$ and $\mu$ in Figures~\ref{batchLim1}. In addition to the interesting shapes of the densities across the different settings, one can see the limiting form of the relationships given by the recursion in Proposition~\ref{recurProp} in these plots. We can note that one could also calculate these through a numerical inverse Laplace transform of the steady-state moment generating function in Corollary~\ref{ssMGF}, although this may likely incur significantly more computational costs than the simulation procedure.

%

So far we have only considered exponentially distributed service. In the next subsection we will address this and extend this Poisson sum representation of the steady-state distribution to hold for general service. We do this through viewing the $n$-batch-size system as being composed of $n$ sub-systems that experience single arrivals simultaneously.

\subsection{Generalizing through Sub-System Perspectives}\label{orderStatSubsec}

Because of the infinite server construction of this model, we can also interpret this system as being a network of sub-systems that also feature infinitely many servers. However, this network's mutuality is not in its services but rather in its arrivals. Specifically, in this subsection we will think of infinite server queues with batch arrivals of size $n$ as being $n$ infinite server queues that all receive individual arrivals simultaneously. From this perspective, one can quickly observe that marginally each subsystem will be distributed as a standard infinite server queue.

For example, if the batch system is the $M_t^n/M/\infty$ that we first considered in Subsection~\ref{transientSub}, then each of these sub-queues are $M_t/M/\infty$ systems. These sub-systems are coupled through the coincidence of their arrival times but otherwise operate independently from one another. To quantify the relationship between these systems, in Proposition~\ref{covgen} we derive the transient covariance between two sub-systems for a general time-varying arrival rate.

\begin{proposition}\label{covgen}
Let the batch arrival queue $Q_t$ with batch size $n \in \mathbb{Z}^+$ be represented as a superposition of $n$ infinite server single arrival queues $\{Q_{t,i} \mid 1 \leq i \leq n \}$ that all receive arrivals simultaneously and each have independent exponentially distributed service, as described above. Let $\lambda(t) > 0$ be the non-stationary rate of simultaneous arrivals and let $\mu > 0$ be the rate of service. Then, for distinct $i , j \in \{1, \dots , n\}$, the covariance between $Q_{t,i}$ and $Q_{t,j}$ is given by
\begin{align}
&
\Cov{Q_{t,i},Q_{t,j}}
=
e^{-2\mu t}
\int_0^t
\lambda(s)
e^{2\mu s}
\mathrm{d}s
\end{align}
for all $t \geq 0$.
\begin{proof}
From Equation~\ref{dynkin}, we can solve for the product moment of the two sub-systems through the ODE
\begin{align*}
\frac{\mathrm{d}}{\mathrm{d}t}\E{Q_{t,i}Q_{t,j}}
&
=
\lambda(t)\left(\E{Q_{t,i}} + \E{Q_{t,j}} + 1\right) - 2\mu \E{Q_{t,i}Q_{t,j}}.
\end{align*}
The solution to this differential equation is given by
\begin{align*}
\E{Q_{t,i}Q_{t,j}}
&
=
Q_{0,i}Q_{0,j}e^{-2\mu t}
+
e^{-2\mu t}
\int_0^t
\lambda(s)
\left(\E{Q_{s,i}}e^{2\mu s} + \E{Q_{s,j}}e^{2\mu s} + e^{2\mu s}\right)
\mathrm{d}s.
\intertext{By substituting the corresponding forms of $\E{Q_{s,k}} = Q_{0,k}e^{-\mu s} + e^{-\mu s} \int_0^s \lambda(u)e^{\mu u}\mathrm{d}u$ in for each of the two means, we have}
\E{Q_{t,i}Q_{t,j}}
&
=
Q_{0,i}Q_{0,j}e^{-2\mu t}
+
e^{-2\mu t}
\int_0^t
\lambda(s)
\bigg(
e^{2\mu s}
+
\left(Q_{0,i} + \int_0^s \lambda(u)e^{\mu u}\mathrm{d}u\right)e^{\mu s}
\\
&
\quad
+
\left(Q_{0,j} + \int_0^s \lambda(u)e^{\mu u}\mathrm{d}u\right)e^{\mu s}
\bigg)
\mathrm{d}s,
\intertext{and this simplifies to the following}
\E{Q_{t,i}Q_{t,j}}
&
=
Q_{0,i}Q_{0,j}e^{-2\mu t}
+
e^{-2\mu t}
\int_0^t
\lambda(s)
e^{2\mu s}
\mathrm{d}s
+
\left(Q_{0,i} + Q_{0,j}\right)
e^{-2\mu t}
\int_0^t
\lambda(s)
e^{\mu s}
\mathrm{d}s
\\
&
\quad
+
2 e^{-2\mu t}
\int_0^t
\lambda(s)
 e^{\mu s}  \int_0^s \lambda(u)e^{\mu u}\mathrm{d}u
\mathrm{d}s.
\intertext{We can now use the fact that for a function $F:\mathbb{R}^+ \to \mathbb{R}$ defined such that $F(t) = \int_0^t f(s) \mathrm{d}s$ for a given $f(\cdot)$, integration by parts implies
$$
\int_0^t f(s) F(s) \mathrm{d}s
=
F(t)^2
-
\int_0^t F(s) f(s)  \mathrm{d}s
,
$$
and so $\int_0^t f(s) F(s) \mathrm{d}s = \frac{F(t)^2}{2}$. This allows us to simplify to}
\E{Q_{t,i}Q_{t,j}}
&
=
Q_{0,i}Q_{0,j}e^{-2\mu t}
+
e^{-2\mu t}
\int_0^t
\lambda(s)
e^{2\mu s}
\mathrm{d}s
+
\left(Q_{0,i} + Q_{0,j}\right)
e^{-2\mu t}
\int_0^t
\lambda(s)
e^{\mu s}
\mathrm{d}s
\\
&
\quad
+
 e^{-2\mu t}
\left(
\int_0^t
\lambda(s)e^{\mu s}
\mathrm{d}s
\right)^2,
\end{align*}
and now we turn our focus to the product of the means. Here we distribute the multiplication to find that
\begin{align*}
\E{Q_{t,i}}\E{Q_{t,j}}
&
=
\left(
Q_{0,i}e^{-\mu t}
+
e^{-\mu t}
\int_0^t
\lambda(s)e^{\mu s}
\mathrm{d}s
\right)
\left(
Q_{0,j}e^{-\mu t}
+
e^{-\mu t}
\int_0^t
\lambda(s)e^{\mu s}
\mathrm{d}s
\right)
\\
&
=
Q_{0,i}Q_{0,j}e^{-2\mu t}
+
(Q_{0,i}+Q_{0,j})e^{-2\mu t}
\int_0^t
\lambda(s)e^{\mu s}
\mathrm{d}s
+
e^{-2\mu t}
\left(
\int_0^t
\lambda(s)e^{\mu s}
\mathrm{d}s
\right)^2
\end{align*}
and by subtracting this expression from that of the product moment, we complete the proof.
\end{proof}
\end{proposition}

As a consequence of this, we can specify the covariance between sub-systems in the non-stationary and stationary arrival settings we have considered thus far in this report. Further, for stationary arrival rates we capitalize on simplified expressions to also give an explicit expression for the correlation coefficient between two sub-systems.

\begin{corollary}\label{covariance}
Let $Q_t$ be an infinite server queue with arrival batch size $n \in \mathbb{Z}^+$ and exponential service rate $\mu > 0$. Further, let $Q_{t,k}$ for $k \in \{1, \dots, n\}$ be infinite server queues with solitary arrivals and exponential service rate $\mu > 0$, so that $\sum_{k=1}^n Q_{t,k} = Q_t$ for all $t \geq 0$. Let $i , j \in \{1, \dots , n\}$ be distinct. Then, if the arrival rate is given by $\lambda + \sum_{k=1}^\infty a_k \cos(k t) + b_k \sin(k t) > 0$, the covariance between $Q_{t,i}$ and $Q_{t,j}$ is
\begin{align}
\Cov{Q_{t,i},Q_{t,j}}
&
=
\frac{\lambda}{2\mu}\left(1 - e^{-2\mu t}\right)
+
\sum_{k=1}^\infty
\frac{a_k}{k^2 + 4\mu^2}
\left(
2\mu \cos(kt) + k \sin(kt) - 2\mu e^{-2\mu t}
\right)
\nonumber
\\
&
\quad
+
\sum_{k=1}^\infty
\frac{b_k}{k^2 + 4\mu^2}\left(2\mu \sin(kt) - k\cos(kt) + k e^{-2\mu t}\right)
,
\end{align}
and if the arrival rate is given by $\lambda > 0$, the covariance between $Q_{t,i}$ and $Q_{t,j}$ is
\begin{align}
&
\Cov{Q_{t,i},Q_{t,j}}
=
\frac{\lambda}{2\mu}\left(1 - e^{-2\mu t}\right)
,
\end{align}
where all $t \geq 0$.  Finally, the correlation between two sub-systems in the stationary setting can be calculated as
\begin{align}
&
\mathrm{Corr}[Q_{t,i},Q_{t,j}]
=
\frac{
\frac{\lambda}{2\mu}\left(1 - e^{-2\mu t}\right)
}
{
\sqrt{
\left(
Q_{0,i}\left(e^{-\mu t} - e^{-2\mu t}\right)
+
 \frac{\lambda }{\mu}
 \left(
 1 - e^{- \mu t}
 \right)
 \right)
 \left(
Q_{0,j}\left(e^{-\mu t} - e^{-2\mu t}\right)
+
 \frac{\lambda }{\mu}
 \left(
 1 - e^{- \mu t}
 \right)
 \right)
 }
}
\nonumber
,
\end{align}
hence for stationary arrival rates, $\mathrm{Corr}[Q_{t,i},Q_{t,j}] \to \frac{1}{2}$ as $t \to \infty$.
\end{corollary}


Thus, we find that for a fully Markovian batch arrival queue with stationary arrival rate the correlation among any two sub-systems in steady-state is $\frac{1}{2}$, regardless of the arrival or service parameters. In some sense this seems to capture a balance between the effect of arrivals and of services on an infinite server system, with the latter being independent between these systems and the former being perfectly correlated.

Now, we can pause to note that we have actually made an implicit modeling choice by separating the batch into $n$ identical sub-systems. In this set-up we have decided to route all customers within one batch equivalently, but we are free to make other routing decisions and still maintain the $n$ sub-systems construction. With that in mind, it seems natural to wonder if we can uncover distributional structure of the full system if we choose our routing procedure carefully. We will now find that not only is this true, but we in fact already have already seen a suggestion on what type of routing to consider.

From Proposition~\ref{ssDist}, we have seen that the steady-state distribution of the $M^n/M/\infty$ system is equivalent to that of $\sum_{j=1}^n j Y_j$ where $Y_j \sim \mathrm{Pois}(\frac{\lambda}{j \mu })$ are independent. We can also note that just as the minimum of the independent sample $S_1, \dots, S_n \sim \mathrm{Exp}(\mu)$ will be exponentially distributed with rate $n\mu$, for $S_{(i)}$ as the $i^\text{th}$ ordered statistic of the $n$-sample we have that $S_{(i)} - S_{(i-1)} \sim \mathrm{Exp}((n-i+1)\mu)$. Of course, the sum of these differences will telescope so that $\sum_{j=1}^i S_{(j)} - S_{(j-1)} = S_{(i)}$.

\begin{figure}[h]
\centering
\vspace{.15in}
\begin{tikzpicture}[line cap=round,line join=round, inner sep = 0pt, outer sep = 0pt]
\draw  (-4,4) rectangle (-3,-2);
\draw  (-1,4) rectangle (0.125,-2);
\node (v1) at (-4,3) {};
\node (v2) at (-3,3) {};
\node at (-3.5,3.5) {$S_{1}$};
\node at (-3.5,2.5) {$S_{2}$};
\node (v4) at (-3,2) {};
\node (v3) at (-4,2) {};
\node (v7) at (-4,-1) {};
\node (v8) at (-3,-1) {};
\node at (-3.5,-1.5) {$S_{n}$};
\node (v5) at (-4,0) {};
\node (v6) at (-3,0) {};
\node at (-3.5,-0.5) {$S_{n-1}$};
\node at (-3.5,1) {\vdots};
\node (v15) at (-1,3) {};
\node (v16) at (.125,3) {};
\node at (-0.42,3.5) {$S_{(1)}$};
\node at (-0.42,2.5) {$S_{(2)}$};
\node (v13) at (-1,2) {};
\node (v14) at (.125,2) {};
\node at (-0.42,1) {\vdots};
\node (v11) at (-1,0) {};
\node (v12) at (.125,0) {};
\node (v9) at (-1,-1) {};
\node (v10) at (.125,-1) {};
\node at (-0.42,-0.5) {$S_{(n-1)}$};
\node at (-0.42,-1.5) {$S_{(n)}$};
\draw  (v1) edge (v2);
\draw  (v3) edge (v4);
\draw  (v5) edge (v6);
\draw  (v7) edge (v8);
\draw  (v9) edge (v10);
\draw  (v11) edge (v12);
\draw  (v13) edge (v14);
\draw  (v15) edge (v16);
\node (v17) at (-6,1) {};
\node (v18) at (-4,1) {};
\node (v19) at (-3,1) {};
\node (v20) at (-1,3.5) {};
\node (v21) at (-1,2.5) {};
\node (v22) at (-1,-0.5) {};
\node (v23) at (-1,-1.5) {};
\draw[-triangle 60]   (v17) edge (v18);
\draw[-triangle 60]   (v19) edge (v20);
\draw[-triangle 60]   (v19) edge (v21);
\draw[-triangle 60]   (v19) edge (v22);
\draw[-triangle 60]   (v19) edge (v23);
\node (v24) at (2.5,3.5) {$Q_1$};
\node (v25) at (2.5,2.5) {$Q_2$};
\node (v26) at (2.5,-0.5) {$Q_{n-1}$};
\node (v27) at (2.5,-1.5) {$Q_{n}$};
\node (v34) at (.125,-1.5) {};
\node (v32) at (.125,-0.5) {};
\node (v30) at (.125,2.5) {};
\node (v28) at (.125,3.5) {};
\draw  (v24) circle (0.5);
\draw  (v25) circle (0.5);
\draw  (v26) circle (0.5);
\draw  (v27) circle (0.5);
\node (v35) at (2,-1.5) {};
\node (v33) at (2,-0.5) {};
\node (v31) at (2,2.5) {};
\node (v29) at (2,3.5) {};
\draw[-triangle 60]   (v28) edge (v29);
\draw[-triangle 60]   (v30) edge (v31);
\draw[-triangle 60]   (v32) edge (v33);
\draw[-triangle 60]   (v34) edge (v35);
\node (v37) at (5,3.5) {};
\node (v39) at (5,2.5) {};
\node (v41) at (5,-0.5) {};
\node (v43) at (5,-1.5) {};
\node (v36) at (3,3.5) {};
\node (v38) at (3,2.5) {};
\node (v40) at (3,-0.5) {};
\node (v42) at (3,-1.5) {};
\draw[-triangle 60]   (v36) edge (v37);
\draw[-triangle 60]   (v38) edge (v39);
\draw[-triangle 60]   (v40) edge (v41);
\draw[-triangle 60]   (v42) edge (v43);
\node at (-5,1.5) {$\lambda$};
\node at (2.5,1) {\vdots};

\node at (-3.5,-2.75) [align=center]{Services\\for Batch};
\node at (-0.5,-2.75) [align=center]{Order\\Statistics};
\node at (2.5,-2.75) [align=center]{Ordered\\Queues};
\end{tikzpicture}
\vspace{.15in}
\caption{Queueing diagram for the batch arrival queue with infinite servers, in which the arriving entities are routed according to the ordering of their service durations.}\label{orderDiagram}
\end{figure}
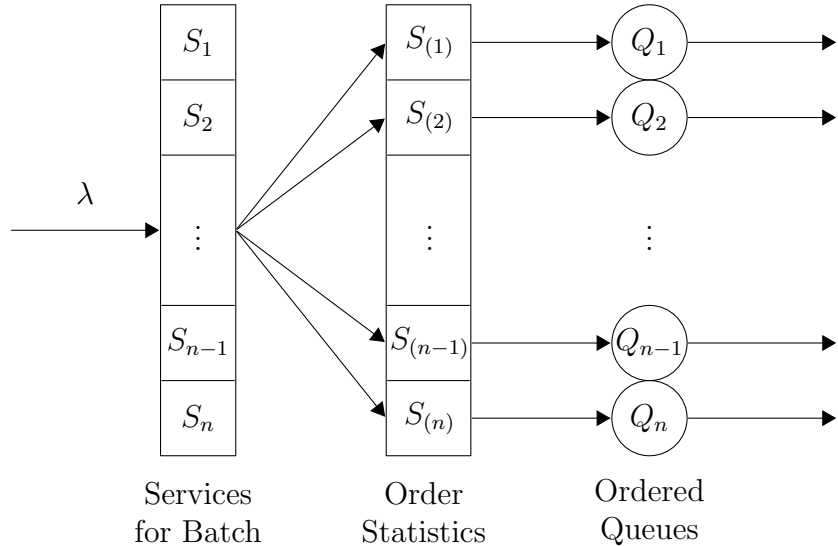

Taking this as inspiration, we will now assume that upon the arrival of a batch we can now know the duration of each customer's service. We then take the sub-queues to be such that the first sub-system always receives the service with the shortest duration, the second sub-system receives the second shortest service, and so on. Thus, we will route each batch of customers according to the order statistics within each batch. For reference, we visualize this sub-system construction in Figure~\ref{orderDiagram}.

We can note that while the covariance structure we explored in Proposition~\ref{covgen} and Corollary~\ref{covariance} do not apply for this new routing, the sub-systems are certainly still correlated. Due to the order-statistics structuring of the service in each queue, we can note that now both the arrival processes and the service distributions will be dependent. However, we can in fact use our understanding of this dependence to not only understand how these systems relate to one another, but also to interpret how they form the structure of the full batch system as a whole. In this way, we will now consider a $M^n/G/\infty$ system. As follows in Theorem~\ref{orderStat}, we will find that the order-statistics-routing inspiration we have used from Proposition~\ref{ssDist} leads us to a generalized Poisson sum result for general service distributions.

\begin{theorem}\label{orderStat}
Let $Q_t(n)$ be an $M^n/G/\infty$ queue. That is, let $Q_t(n)$ be an infinite server queue with stationary arrival rate $\lambda > 0$, arrival batch size $n \in \mathbb{Z}^+$, and general service distribution $G$. Then, the steady-state distribution of the number in system $Q_\infty(n)$ is
\begin{align}
Q_\infty(n)
\stackrel{D}{=}
\sum_{j=1}^n
(n-j+1) Y_j
\end{align}
where $Y_j \sim \mathrm{Pois}\left(\lambda \E{S_{(j)} - S_{(j-1)}}\right)$ are independent, with $S_{(1)} \leq \dots \leq S_{(n)}$ as order statistics of the distribution $G$ and with $S_{(0)} = 0$.
\begin{proof}
As we have discussed in the paragraphs preceding this statement, we will consider the full queueing system as being composed of $n$ infinite server sub-systems to which we route the arriving customers in each batch. That is, let $Q_1$, \dots, $Q_n$ be infinite server queues of which we will consider the steady-state behavior. Upon the arrival of a batch, we order the customers according to the duration of their service.  Then, we send the customer with the earliest service completion to $Q_1$, the customer with the second earliest to $Q_2$, and so on.

When viewing each sub-system on its own, we see that $Q_j$ is an infinite server queue with single arrivals according to a Poisson process with rate $\lambda$ and service distribution matching that of $S_{(j)}$, the $j^\text{th}$ order statistics of $G$. Thus, we can see that in steady-state $Q_j \sim \mathrm{Pois}\left(\lambda \E{S_{(j)}}\right)$ through the literature for $M/G/\infty$ queues, such as in \cite{eick1993physics}. While we can further observe that $Q_\infty(n) = \sum_{j=1}^n Q_j$, we must take care in re-assembling the sub-queues. In particular, we can note that $S_{(j)}$ shares a similar structure with $S_{(j-1)}$. Each order statistic can be viewed as a construction of the gaps between the lower ordered quantities:
$$
S_{(j)} = \sum_{k=1}^j S_{(k)} - S_{(k-1)} .
$$
Thus, from the thinning property of the Poisson distribution and the linearity of expectation, we can write the distribution of $Q_j$ as a sum of independent Poisson RV's, as given by
$$
Q_j
\sim
 \sum_{k=1}^j
\mathrm{Pois}\left(\lambda \E{S_{(k)} - S_{(k-1)}}\right)
.
$$
We can note further that $j-1$ of the Poisson components of $Q_j$ are the exact components of $Q_{j-1}$, with $j-2$ of these components also shared with $Q_{j-2}$, $j-3$ with $Q_{j-3}$, and so on. Then, we see that the Poisson component $\mathrm{Pois}\left(\lambda \E{S_{(j)} - S_{(j-1)}}\right)$ is repeated $n-j+1$ times across this sub-system construction of $Q_\infty(n)$, as it appears in each of the Poisson sum expressions of $Q_{j}$, $Q_{j+1}$, \dots, $Q_{n-1}$, and $Q_n$. Assembling $Q_\infty(n)$ in this way, we complete the proof.
\end{proof}
\end{theorem}

One can also note that this order statistic sub-system structure also provides some motivation for the occurrence of the harmonic numbers that we observed in Subsection~\ref{markovSS} when viewing the largest order statistic, which we discuss now in the following remark.
\begin{remark}
For $S_i \sim \mathrm{Exp}(\mu)$, one can see through the telescoping construction of the order statistics that
$$
\E{S_{(n)}}
=
\sum_{i=1}^n \E{S_{(i)} - S_{(i-1)}}
=
\sum_{i=1}^n \frac{1}{(n-i+1) \mu}
=
\frac{1}{\mu} H_n
.
$$
\end{remark}

Now, throughout this section we have operated on the assumption that the batch size is a known, fixed constant. While this may be applicable in some settings there are certainly many settings where the batch size is unknown and varies between arrivals. Thus, we address this in Section~\ref{randomSec} and find that many of the results we have shown thus far can be replicated for models with random batch size.

\section{Random Batch Sizes}\label{randomSec}

We will now consider systems in which the size of an arriving batch is drawn from an independent and identically distributed sequence of random variables. We will treat the distribution of the batch size as general throughout this work. As in Section~\ref{deterSec}, we assume that the times of arrivals are given by a Poisson process, with consideration given to both stationary and non-stationary rates, and we will again analyze both exponential and general service distributions.

We start by giving the mean and variance of the system for time-varying arrival rates with exponential service in Subsection~\ref{meanvarRandSubsec}. Then, in Subsection~\ref{limitRandSubsec} we give three limiting results for the stationary arrivals model: a batch scaling, a fluid limit, and a diffusion limit. Finally in Subsection~\ref{orderStatRandSubsec} we extend the Poisson sum construction of the steady-state distribution to hold for random batch sizes.

One can note that many of these results are generalizations or extensions of findings from Section~\ref{deterSec}, thus implying them as a special case and perhaps even building a case for them to be omitted. Rather, these findings are critical to the narrative of this report. As we will see, the results for fixed batch size provide the analytic foundations and conceptual inspirations from which we derive much of the analysis in this section.

\subsection{Mean and Variance for Time-Varying, Markovian Case}\label{meanvarRandSubsec}

To begin our exploration into random batch size systems, we'll start simple: we'll look at a fully Markovian (albeit time-varying) system and find the mean and variance, using conditional probability and our results from Section~\ref{deterSec}. Specifically, in this subsection we will consider the $M_t^N/M/\infty$ queue. That is, take an infinite server queue with a general non-stationary arrival rate. We suppose that arrivals occur in batches of random size from a sequence of independent and identically distributed random variables. Furthermore, we suppose that service is exponentially distributed. We now give the mean and variance of this system in Proposition~\ref{meanvarRand}.

\begin{proposition}\label{meanvarRand}
Let $Q_t$ be an infinite server queue with finite, time-varying arrival rate $\lambda(t) > 0$, exponential service rate $\mu > 0$, and random batch size with finite mean, $\E{N}$. Then, the mean number in system is given by
\begin{align}
\E{Q_t}
&
=
Q_0 e^{-\mu t}
+
e^{-\mu t}\E{N}
\int_0^t \lambda(s)e^{\mu s}\mathrm{d}s
 ,
\end{align}
for all $t \geq 0$. Then, if the batch size distribution has finite second moment $\E{N^2}$, the variance of the number in system is given by
\begin{align}
\Var{Q_t}
&
=
Q_0\left(e^{-\mu t} - e^{-2\mu t}\right)
+
e^{-2\mu t}\left( \E{N^2} - \E{N} \right)
\int_0^t \lambda(s) e^{2\mu s} \mathrm{d}s
\nonumber
\\
&
\quad
+
e^{-\mu t}\E{N}\int_0^t \lambda(s) e^{\mu s} \mathrm{d}s
,
\end{align}
again for all $t \geq 0$.
\begin{proof}
Using the infinitesimal generator method, we have that the first and second moments of this system are given by the solutions to
\begin{align*}
\frac{\mathrm{d}}{\mathrm{d}t}\E{Q_t}
&
=
\lambda(t) \E{N_1} - \mu \E{Q_t}
,
\\
\frac{\mathrm{d}}{\mathrm{d}t}\E{Q_t^2}
&
=
\lambda(t) \left(2\E{Q_t}\E{ N_1}+\E{N_1^2}\right) - 2\mu \E{Q_t^2} + \mu \E{Q_t}
,
\end{align*}
where $\{N_i \mid i \in \mathbb{Z}^+\}$ are the i.i.d. batch sizes that are also independent of the queue. Through noting that
$$
\frac{\mathrm{d}}{\mathrm{d}t}\Var{Q_t}
=
\frac{\mathrm{d}}{\mathrm{d}t}\E{Q_t^2}
-
2\E{Q_t} \frac{\mathrm{d}}{\mathrm{d}t}\E{Q_t}
=
\lambda(t)\E{N_1^2} + \mu \E{Q_t} - 2\mu \Var{Q_t}
,
$$
we can solve for the stated results.
\end{proof}
\end{proposition}

In addition to providing a direct comparison to the fixed batch size case in conjunction with Corollary~\ref{meanvarcor}, Proposition~\ref{meanvarRand} also provides a building block for the remainder of this section. In particular, in the following subsection we will develop a series of limiting results for this queueing system, including fluid and diffusion limits. In those cases, we will use this result for added interpretation. To expedite comparison in cases of stationary arrival rates, we now give the mean and variance for such systems in Corollary~\ref{meanvarRandStat}. Additionally, to also facilitate comparison to Corollary~\ref{meanvarcor}, we provide expressions for periodic arrival rates in Corollary~\ref{meanvarRandPer}.

\begin{corollary}\label{meanvarRandStat}
Let $Q_t$ be an infinite server queue with stationary arrival rate $\lambda > 0$, exponential service rate $\mu > 0$, and random batch size with mean $\E{N}$. Then, the mean number in system is given by
\begin{align}
\E{Q_t}
&
=
Q_0 e^{-\mu t}
+
 \frac{\lambda \E{N}}{\mu}
 \left(
 1 - e^{- \mu t}
 \right)
 ,
\end{align}
for all $t \geq 0$. Then, if the batch size distribution has finite second moment $\E{N^2}$, the variance of the number in system is given by
\begin{align}
\Var{Q_t}
&
=
Q_0\left(e^{-\mu t} - e^{-2\mu t}\right)
+
 \frac{\lambda \E{N}}{\mu}
 \left(
 1 - e^{- \mu t}
 \right)
 +
 \frac{ \lambda }{2\mu}
  \left(
 \E{N^2} - \E{N}
 \right)
 \left(
 1 - e^{-2 \mu t}
 \right)
,
\end{align}
again for all $t \geq 0$.
\end{corollary}

\begin{corollary}\label{meanvarRandPer}
Let $Q_t$ be an infinite server queue with periodic arrival rate $\lambda + \sum_{k=1}^\infty a_k \cos(k t) + b_k \sin(k t) > 0$, exponential service rate $\mu > 0$, and random batch size with finite mean, $\E{N}$. Then, the mean number in system is given by
\begin{align}
\E{Q_t}
&
=
Q_0 e^{-\mu t}
+
 \frac{\lambda \E{N}}{\mu}
 \left(
 1 - e^{- \mu t}
 \right)
 +
 \sum_{k=1}^\infty
 \frac{\E{N}(a_k  \mu - b_k k)  }{k^2+\mu^2}
  \left(
 \cos(kt) - e^{-\mu t}
 \right)
 \nonumber
 \\
 &
 \quad
 +
 \sum_{k=1}^\infty
 \frac{\E{N}(a_k k  +  b_k  \mu )}{k^2+\mu^2} \sin(k t)
 ,
\end{align}
for all $t \geq 0$. Then, if the batch size distribution has finite second moment $\E{N^2}$, the variance of the number in system is given by
\begin{align}
\Var{Q_t}
&
=
Q_0\left(e^{-\mu t} - e^{-2\mu t}\right)
+
 \frac{\lambda \E{N}}{\mu}
 \left(
 1 - e^{- \mu t}
 \right)
 +
 \sum_{k=1}^\infty
 \frac{E{N}(a_k  \mu - b_k k)  }{k^2+\mu^2}
  \left(
 \cos(kt) - e^{-\mu t}
 \right)
 \nonumber
 \\
 &
 \quad
 +
 \sum_{k=1}^\infty
 \frac{\E{N}(a_k k  +  b_k  \mu )}{k^2+\mu^2} \sin(k t)
 +
 \frac{ \lambda }{2\mu}
  \left(
 \E{N^2} - \E{N}
 \right)
 \left(
 1 - e^{-2 \mu t}
 \right)
 \nonumber
 \\
 &
 \quad
 +
  \left(
 \E{N^2} - \E{N}
 \right)
 \left(
 \sum_{k=1}^\infty
 \frac{2 a_k \mu - b_k k  }{k^2 + 4\mu^2}
  \left(
 \cos(kt) - e^{-2\mu t}
 \right)
 +
 \sum_{k=1}^\infty
 \frac{a_k k  +  2 b_k \mu }{k^2+ 4\mu^2} \sin(k t)
 \right)
,
\end{align}
again for all $t \geq 0$.
\end{corollary}

\subsection{Limiting Results for Stationary Arrival Rates}\label{limitRandSubsec}

We will now focus on systems with stationary arrival rates throughout the analysis in this subsection. In doing so, we derive limit theorems for various scalings of this process. To begin, we show a brief technical lemma for the limit of non-negative random variables that can be represented as sums of independent and identically distributed random variables.

\begin{lemma}\label{divisLemma}
Let $X(n)$ be any random variable that
$
X(n) = \sum_{k=1}^n Y_k
$
where $Y_k$ are i.i.d. non-negative, discrete random variables. Then, the moment generating function of $X(n)$ is such that
\begin{align*}
\E{e^{\frac{\theta X(n)}{n}}}
\to
e^{\E{Y_1} \theta}
\end{align*}
as $n \to \infty$.
\begin{proof}
By the strong law of large numbers, we have that
$$
\lim_{n\to\infty}\frac{X(n)}n = \lim_{n\to\infty} \frac{1}{n}\sum_{k=1}^n Y_k \stackrel{\text{a.s.}}{=} \E{Y_1},
$$
and this implies convergence in distribution, which is equivalent to convergence of moment generating functions.
\end{proof}
\end{lemma}

We can note that this condition is a weaker form of infinite divisibility. Thus, in addition to holding for any infinitely divisible and non-negative random variables such as the Poisson, and negative binomial distributions, Lemma~\ref{divisLemma} also holds for some distributions that are not infinitely divisible, such as the binomial. Using this lemma we can now find our first limit theorem for random batch sizes, a batch scaling result akin to Proposition~\ref{scaleMGF}.

\begin{theorem}\label{scaleMGFrand}
For $n \in \mathbb{Z}^+$, let $Q_t(n)$ be an infinite server queue with batch arrivals where the batch size is drawn from the i.i.d. sequence $\{N_i(n) \mid i \in \mathbb{Z}^+\}$. Let $\lambda > 0$ be the arrival rate and let $\mu > 0$ be the rate of exponentially distributed service. Then, suppose that for any $i$ and $n$ there is a sequence of i.i.d. non-negative, discrete random variables $\{B_k \mid k \in \mathbb{Z}^+\}$ such that
$
N_i(n)
=
\sum_{k=1}^n
B_k
.
$
Then, the limiting moment generating function of the batch scaled object
\begin{align}
\lim_{n\to\infty}
\E{e^{\frac{\theta}{n}Q_t(n)}}
=
\begin{cases}
e^{
\frac{\lambda}{\mu}
 \left(
\mathrm{Ei}\left(\theta\E{B_1}\right)
-
\mathrm{Ei}\left(\theta\E{B_1} e^{-\mu t} \right)
-
\mu t
 \right)
}
&
\text{ if $\theta > 0$,}
\\
e^{
\frac{\lambda}{\mu}
 \left(
E_1\left(-\theta\E{B_1} e^{-\mu t}\right)
-
E_1\left(-\theta\E{B_1} \right)
-
\mu t
 \right)
}
&
\text{ if $\theta < 0$,}
\\
1
&
\text{ if $\theta = 0$,}
\end{cases}
\end{align}
for all $t \geq 0$.
\begin{proof}
Because this system is Markovian, we can calculate the time derivative of the moment generating function for a given $n$ as
\begin{align*}
\frac{\mathrm{d}}{\mathrm{d}t}\E{e^{\frac{\theta}{n} Q_t(n)}}
&
=
\E{
\lambda \left(e^{\frac{\theta}{n}N_1(n)}-1\right) e^{\frac{\theta}{n} Q_t(n)}
+
\mu Q_t(n) \left(e^{-\frac{\theta}{n}}-1\right) e^{\frac{\theta}{n} Q_t(n)}
}
\\
&
=
\lambda \left(\E{e^{\frac{\theta}{n}N_1(n)}}-1\right) \E{e^{\frac{\theta}{n} Q_t(n)}}
+
n \mu \left(e^{-\frac{\theta}{n}}-1\right) \E{\frac{Q_t(n)}{n}e^{\frac{\theta}{n} Q_t(n)}}
.
\intertext{This can then be re-expressed in partial differential equation form as}
\frac{\partial \mathcal{M}^n(\theta, t)}{\partial t}
&
=
\lambda \left(\E{e^{\frac{\theta}{n}N_1(n)}}-1\right)  \mathcal{M}^n(\theta, t)
+
n \mu \left(e^{-\frac{\theta}{n}}-1\right) \frac{\partial \mathcal{M}^n(\theta, t)}{\partial \theta},
\end{align*}
where $\mathcal{M}^n(\theta, t) = \E{e^{\frac{\theta}n Q_t(n)}}$. Now, through Lemma~\ref{divisLemma}, we see that the limit of this partial differential equation is given by
\begin{align*}
\frac{\partial \mathcal{M}^\infty(\theta, t)}{\partial t}
&
=
\lambda \left(e^{\theta \E{B_1}}-1\right)  \mathcal{M}^\infty(\theta, t)
-
 \mu \theta \frac{\partial \mathcal{M}^\infty(\theta, t)}{\partial \theta}
 .
\end{align*}
We achieve the stated result through a straightforward update of the method of characteristics approach in Proposition~\ref{scaleMGF}.
\end{proof}
\end{theorem}

We can note that a similar batch scaling of infinite server queues is discussed in \citet{de2017shot}, in which the authors show that the limiting process can be interpreted as a shot noise process. However, that work considers a different class of batch size distributions, as the authors define their batch size distribution in terms of the distribution of the marks through use of a ceiling rounding function. In this way, that paper is more oriented around the distribution of the marks in the shot noise process rather than the size of the batches.

From this result, we can identify a relationship between the moment generating functions of the deterministic and random batch size queues under batch scalings. Let $\mathcal{M}_n^\infty(\theta, t)$ be the limiting moment generating function of the fixed batch size queue as given in Proposition~\ref{scaleMGF} and let $\mathcal{M}_N^\infty(\theta,t)$ be the same for the random batch size queue as we have now seen in Theorem~\ref{scaleMGFrand}. Then, we can observe that
$$
\mathcal{M}_N^\infty(\theta,t)
=
\mathcal{M}_n^\infty(\theta \E{B_1},t)
,
$$
whenever the distribution of the random batch sizes meets the ``finite divisibility'' condition as described in Lemma~\ref{divisLemma}. The relationship between these limiting objects provides a direct comparison between the two different batch types.

As two additional limiting results, we now provide fluid and diffusion limits for scaling the arrival rate in Theorems~\ref{batchFluid} and~\ref{batchDiff}, respectively. We did not give fluid or diffusion limits for the deterministic batch cases in Section~\ref{deterSec}, so these two limits are built from scratch within this section. Although we did not develop such limits explicitly for the $M^n/M/\infty$ system, we will find that these limits can still be used to draw comparisons between this system and the $M^N/M/\infty$  queue simply by treating the random batch size as deterministically distributed. We now begin with the fluid limit.

\begin{theorem}\label{batchFluid}
For $n \in \mathbb{Z}^+$, let $Q_t(n)$ be an infinite server queue with batch arrivals where the batch size is drawn from the i.i.d. sequence $\{N_i \mid i \in \mathbb{Z}^+\}$. Let $n \lambda > 0$ be the arrival rate and let $\mu > 0$ be the rate of exponentially distributed service. Then, the limiting moment generating function of the fluid scaling is given by
\begin{align}
\lim_{n \to \infty} \E{e^{\frac{\theta}{n}Q_t(n)}}
&
=
e^{ \frac{\lambda \E{N_1} \theta}{\mu} \left(1 - e^{-\mu t}\right) + Q_0 \theta  e^{-\mu t}}
,
\end{align}
for all $t \geq 0$.
\begin{proof}
We begin with the infinitesimal generator equation for the time derivative of the moment generating function at a given $n$. This is
\begin{align*}
\frac{\mathrm{d}}{\mathrm{d}t}\E{e^{\frac{\theta}{n}Q_t(n)}}
&
=
\E{
n \lambda \left(e^{\frac{\theta N_1}{n}} - 1\right) e^{\frac{\theta}{n} Q_t(n)}
+
\mu Q_t(n) \left(e^{-\frac{\theta}{n}} - 1\right) e^{\frac{\theta}{n} Q_t(n) }
}
\\
&
=
n \lambda \left(\E{e^{\frac{\theta N_1}{n}}} - 1\right) \E{e^{\frac{\theta}{n} Q_t(n)}}
+
\mu n \left(e^{-\frac{\theta}{n}} - 1\right) \E{\frac{Q_t(n)}{n} e^{\frac{\theta}{n} Q_t(n) }}
,
\intertext{which can also be expressed in partial differential equation form as}
\frac{\partial\mathcal{M}^n(\theta,t)}{\partial t}
&
=
n \lambda \left(\E{e^{\frac{\theta N_1}{n}}} - 1\right) \mathcal{M}^n(\theta, t)
+
\mu n \left(e^{-\frac{\theta}{n}} - 1\right) \frac{\partial\mathcal{M}^n(\theta,t)}{\partial \theta}
,
\end{align*}
where $M^n(\theta,t) = \E{e^{\frac{\theta}n Q_t(n)}}$. By a Taylor expansion of the function $e^{\frac{\theta N_1}n}$ and by taking the limit as $n \to \infty$, we can see that this yields
\begin{align*}
\frac{\partial\mathcal{M}^\infty(\theta,t)}{\partial t}
&
=
\lambda \theta \E{N_1}  \mathcal{M}^\infty(\theta, t)
-
\mu \theta \frac{\partial\mathcal{M}^\infty(\theta,t)}{\partial \theta}
.
\end{align*}
Using the initial condition $\mathcal{M}^\infty(\theta,0) = e^{Q_0 \theta}$, we can see that the solution to this partial differential equation will be
$$
\mathcal{M}^\infty(\theta,t)
=
e^{ \frac{\lambda \E{N_1} \theta}{\mu} \left(1 - e^{-\mu t}\right) + Q_0 \theta  e^{-\mu t}}
,
$$
and this completes the proof.
\end{proof}
\end{theorem}

From Corollary~\ref{meanvarRandStat}, we see that the mean number in system for the $M^N/M/\infty$ queue is $\frac{\lambda \E{N_1} }{\mu} \left(1 - e^{-\mu t}\right) + Q_0  e^{-\mu t}$. Thus, this fluid limit moment generating function is equivalent to $e^{\theta \E{Q_t}}$ for all $t \geq 0$ and all $\theta$, showing that the fluid limit converges to the mean. We now find a connection to both the mean and the variance through a diffusion limit in Theorem~\ref{batchDiff}.

\begin{theorem}\label{batchDiff}
For $n \in \mathbb{Z}^+$, let $Q_t(n)$ be an infinite server queue with batch arrivals where the batch size is drawn from the i.i.d. sequence $\{N_i \mid i \in \mathbb{Z}^+\}$. Let $n \lambda > 0$ be the arrival rate and let $\mu > 0$ be the rate of exponentially distributed service. Then, the limiting moment generating function of the diffusion scaling is given by
\begin{align}
\lim_{n \to \infty} \E{e^{\frac{\theta}{\sqrt{n}}\left(Q_t(n) - \frac{n\lambda \E{N_1}}{\mu}\right)}}
&
=
e^{\frac{\lambda \theta^2}{4\mu} \left(\E{N_1}+\E{N_1^2}\right) \left(1-e^{-\mu t}\right) + \theta Q_0 e^{-\mu t}}
\end{align}
which gives a steady-state approximation of $X \sim \mathrm{Norm}\left(\frac{\lambda \E{N_1}}{\mu}, \frac{\lambda }{2\mu} \left(\E{N_1}+\E{N_1^2}\right) \right)$.
\begin{proof}
Through use of the infinitesimal generator, we have that the time derivative of the moment generating function for a given $n$ can be expressed
\begin{align*}
&
\frac{\mathrm{d}}{\mathrm{d}t}\E{e^{\frac{\theta}{\sqrt{n}}\left(Q_t(n) - \frac{n\lambda \E{N_1}}{\mu}\right)}}
\\
&
=
\E{
n \lambda \left(e^{\frac{\theta N_1}{\sqrt{n}}} - 1\right) e^{\frac{\theta}{\sqrt{n}}\left(Q_t(n) - \frac{n\lambda \E{N_1}}{\mu}\right)}
+
\mu Q_t(n) \left(e^{-\frac{\theta}{\sqrt{n}}} - 1\right) e^{\frac{\theta}{\sqrt{n}}\left(Q_t(n) - \frac{n\lambda \E{N_1}}{\mu}\right)}
}
\\
&
=
\E{\sqrt{n} \lambda \left( \theta N_1 + \frac{\theta^2 N_1^2}{2\sqrt{n}} + \mathrm{O}\left( \frac{\theta^3 N_1^3}{6 n}\right) \right) e^{\frac{\theta}{\sqrt{n}}\left(Q_t(n) - \frac{n\lambda \E{N_1}}{\mu}\right)}}
\\
&
\quad
+
\E{\mu \sqrt{n} \left(\frac{Q_t(n)}{\sqrt{n}} - \frac{n \lambda \E{N_1}}{\sqrt{n} \mu} + \frac{n \lambda \E{N_1}}{\sqrt{n} \mu} \right) \left(e^{-\frac{\theta}{\sqrt{n}}} - 1\right) e^{\frac{\theta}{\sqrt{n}}\left(Q_t(n) - \frac{n\lambda \E{N_1}}{\mu}\right)}}
,
\end{align*}
where here we have used a Taylor expansion of the function $e^{\frac{\theta N_1}{\sqrt{n}}}$. Now, for $\mathcal{M}^n(\theta, t) = \E{e^{\frac{\theta}{\sqrt{n}}\left(Q_t(n) - \frac{n\lambda \E{N_1}}{\mu}\right)}}$, this equation can be written as a partial differential equation as follows:
\begin{align*}
\frac{\partial \mathcal{M}^n(\theta, t)}{\partial t}
&
=
\lambda \theta \sqrt{n} \E{N_1} \mathcal{M}^n(\theta, t)
+
\frac{\lambda \theta^2}{2} \E{N_1^2}\mathcal{M}^n(\theta, t)
+
\sqrt{n}\lambda \E{\mathrm{O}\left(\frac{\theta^3 N_1^3}{6 n}\right) e^{\frac{\theta}{\sqrt{n}}\left(Q_t(n) - \frac{n\lambda \E{N_1}}{\mu}\right)}}
\\
&
\quad
+
\sqrt{n}\mu \left(e^{-\frac{\theta}{\sqrt{n}}}-1\right) \frac{\partial \mathcal{M}^n(\theta, t)}{\partial \theta}
+
n \lambda \E{N_1} \left(e^{-\frac{\theta}{\sqrt{n}}}-1\right) \mathcal{M}^n(\theta, t)
.
\end{align*}
As we take $n \to \infty$ this PDE becomes
\begin{align*}
\frac{\partial \mathcal{M}^\infty(\theta, t)}{\partial t}
&
=
\frac{\lambda \theta^2}2  \E{N_1} \mathcal{M}^\infty(\theta, t)
+
\frac{\lambda \theta^2}{2} \E{N_1^2}\mathcal{M}^\infty(\theta, t)
-
\mu \theta \frac{\partial \mathcal{M}^\infty(\theta, t)}{\partial \theta}
,
\end{align*}
and this yields a solution of
$$
\mathcal{M}^\infty(\theta, t)
=
e^{\frac{\lambda \theta^2}{4\mu} \left(\E{N_1}+\E{N_1^2}\right) \left(1-e^{-\mu t}\right) + \theta Q_0 e^{-\mu t}}
.
$$
To observe the steady-state distribution, we take the limit as $t \to \infty$ and observe that this produces the moment generating function for a Gaussian.
\end{proof}
\end{theorem}

By comparison to the limits of the expresions in Corollary~\ref{meanvarRandStat} as $t \to \infty$, we can now observe that this steady-state approximation is equal in mean and variance to the steady-state queue.


\subsection{Extending the Order Statistics Sub-Systems}\label{orderStatRandSubsec}

In Subsection~\ref{orderStatSubsec} we found that the steady-state distribution of infinite server queues with fixed batch size and general service can be written as a sum of scaled Poisson random variables, providing a succinct interpretation of the process and an efficient simulation procedure for approximate calculations. The underlying observation that supported this approach was that we can think of an infinite server queue with batch arrivals as a collection of infinite server queues with solitary arrivals that occur simultaneously. Using the thinning property of Poisson processes, we now extend this result to queues with random batch sizes and general service.

\begin{theorem}\label{orderStatRand}
Let $Q_t$ be a $M^N/G/\infty$ queue. That is, let $Q_t$ an infinite server queue with stationary arrival rate $\lambda > 0$, arrival batch of random size according to the i.i.d.~sequence of non-negative integer valued random variables $\{N_i \mid i \in \mathbb{Z}^+\}$, and general service distribution $G$. Then, the steady-state distribution of the number in system $Q_\infty$ is
\begin{align}
Q_\infty
\stackrel{D}{=}
\sum_{n=1}^\infty \sum_{j=1}^n
(n-j+1) Y_{j,n}
\end{align}
where $Y_{j,n} \sim \mathrm{Pois}\left(\lambda p_n \E{S_{(j,n)} - S_{(j-1,n)}}\right)$ are independent, with $S_{(1,n)} \leq \dots \leq S_{(n,n)}$ as order statistics of the distribution $G$ when $N_i = n$, where $S_{(0,n)} = 0$ for all $n$ and $p_n = \PP{N_1 = n}$.
\begin{proof}
To begin, we suppose that there is some $m \in \mathbb{Z}^+$ such that $\PP{N_i \in \{0,\dots, m\}} = 1$. Then, using the thinning property of Poisson processes, we separate the arrival process into $m$ arrival streams where the $n^\text{th}$ arrival rate is $\lambda p_n$. Then, by Theorem~\ref{orderStat} the steady-state distribution of the number in system from the $n^\text{th}$ stream is
$$
\sum_{j=1}^n
(n-j+1) \mathrm{Pois}\left(\lambda p_n \E{S_{(j,n)} - S_{(j-1,n)}}\right)
.
$$
Then, since the $m$ thinned Poisson streams are independent, we have that the full combined system will be distributed as
$$
\sum_{n=1}^m
\sum_{j=1}^n
(n-j+1) \mathrm{Pois}\left(\lambda p_n \E{S_{(j,n)} - S_{(j-1,n)}}\right)
.
$$
Through taking the limit as $m \to \infty$, we achieve the stated result.
\end{proof}
\end{theorem}

We can note that Theorem~\ref{orderStatRand} also provides a method for approximate empirical calculation through simulation. This representation can also be simplified if more information is known about the distribution of the batch size or of the service, or both. As an example, we give the distribution for the fully Markovian system in the following corollary.

\begin{corollary}
Let $Q_t$ be a $M^N/M/\infty$ queue. That is, let $Q_t$ an infinite server queue with stationary arrival rate $\lambda > 0$, arrival batch of random size according to the i.i.d.~sequence of non-negative integer valued random variables $\{N_i \mid i \in \mathbb{Z}^+\}$, and exponentially distributed service at rate $\mu > 0$. Then, the steady-state distribution of the number in system $Q_\infty$ is
\begin{align}
Q_\infty
\stackrel{D}{=}
\sum_{j=1}^\infty
j Y_{j}
\end{align}
where $Y_{j} \sim \mathrm{Pois}\left(\frac{\lambda}{j \mu} \bar F_N(j) \right)$ are independent, where $\bar F_N(j) = \PP{N_1 \geq j}$.
\end{corollary}

One can note that the moment generating function for this system in steady-state is
$$
\E{e^{\theta Q_\infty}}
=
e^{\sum_{j=1}^\infty \frac{\lambda}{j\mu} \bar{F}_N(j)\left(e^{j \theta} - 1\right)}
,
$$
and that this also admits a connection to the generalized Hermite distributions we discussed in Subsection~\ref{markovSS}. In particular, this generalized Hermite distribution can be characterized by $\frac{\lambda}{\mu}$, which is again the mean of the distribution, and the complementary cumulative distribution function of the batch size distribution, which dictates the coefficients at each $j$. For this reason, it may be possible that the steady-state distribution of the queue may be simplified even further for particular batch size distributions.

Because Theorem~\ref{orderStatRand} is again built upon an order statistics sub-queue perspective, it is natural to wonder how the distribution of the batch size would affect those sub-systems. In particular, we now consider the following scenario: suppose that the batch size is bounded by some constant, say $k$, and that we have $k$ sub-systems. For each arriving batch, the customer with the shortest service duration will go to the first sub-system, the second shortest to the second sub-system, and so on, but only up to the number that have just arrived: if this batch is of size $k-1$, the $k^\text{th}$ sub-queue will not receive an arrival. In this way, the $i^\text{th}$ sub-queue represents the number in system that were the $i^\text{th}$ smallest in their batch. In the following proposition we find the conditions on the batch size distribution under which the distributions of the sub-queues will be equivalent.

\begin{proposition}
Consider a $M^B/G/\infty$ queueing system in which the distribution of $B$ has support on $\{1,\dots, k\}$. Let $\phi \in [0,1]^{k-1}$ be such that $\phi_i = \PP{B = i}$, yielding $\PP{B = k} = 1 - \sum_{i=1}^{k-1} \phi_i$. Let $S_{(i,j)}$ be the $i^\text{th}$ order statistics in a sample of size $j$ from the service distribution.  Furthermore, let $Q_i$ be steady-state number in system of an infinite server sub-queue to which the customer with the $i^\text{th}$ smallest service duration in an arriving batch will be routed whenever there are at least $i$ customers in the batch. Let $M \in \mathbb{R}^{k -1 \times k-1}$ be an upper triangular matrix such that
$$
M_{i,j} = \frac{\E{S_{(i,j)}}}{\E{S_{(k,k)}} - \E{S_{(i,k)}}}
,
$$
for $i \leq j$ and $M_{i,j} = 0$ otherwise. For $\mathbf{v} \in \mathbb{R}^{k-1}$ as the all-ones column vector, if $\phi$ is such that
$$
\mathbf{v}
=
\left(M + \mathbf{v} \mathbf{v}^\T \right) \phi
,
$$
then $Q_i \stackrel{D}{=} Q_j$ for all sub-queues $i$ and $j$. Moreover, if $1 + \mathbf{v}^\T M^{-1}\mathbf{v} \ne 0$, then the distributions of the sub-queues are equivalent if and only if $\phi = (M + \mathbf{v}\mathbf{v}^\T)^{-1}\mathbf{v}$.
\begin{proof}
We start by considering the mean of each queue and solving for $\phi$ such that all the means are equal. Let $\lambda$ be the batch arrival rate. Then, the mean of $Q_i$ is
\begin{align*}
\E{Q_i}
&
=
\sum_{j=i}^{k-1} \lambda \phi_j \E{S_{(i,j)}}
+
\lambda \left(1 - \sum_{j=1}^{k-1}\phi_j\right) \E{S_{(i,k)}}
,
\end{align*}
as entities only arrive to $Q_i$ when $B \geq i$. We can note that for $Q_k$ this is
$$
\E{Q_k}
=
\lambda \left(1 - \sum_{j=1}^{k-1}\phi_j\right) \E{S_{(k,k)}}
.
$$
Then, we can see that all the queue means will be equal if $\E{Q_i} = \E{Q_k}$ for all $i$. Thus, we want to solve for $\phi$ such that
$$
0
=
\sum_{j=i}^{k-1} \lambda \phi_j \E{S_{(i,j)}}
+
\lambda \left(1 - \sum_{j=1}^{k-1}\phi_j\right) \E{S_{(i,k)}}
-
\lambda \left(1 - \sum_{j=1}^{k-1}\phi_j\right) \E{S_{(k,k)}}
,
$$
for all $i$. Rearranging this equation and dividing by $\lambda(\E{S_{(k,k)}} - \E{S_{(i,k)}})$, we receive
$$
\sum_{j=i}^{k-1}
\frac{\E{S_{(i,j)}}}{\E{S_{(k,k)}} - \E{S_{(i,k)}}} \phi_j
+
\sum_{j=1}^{k-1}
\phi_j
=
1
.
$$
We can now observe that this forms the linear system $(M + \mathbf{v}\mathbf{v}^\T)\phi = \mathbf{v}$, and so we have shown that if $\phi$ satisfies this system then the means of the sub-queues will be equal. We can note moreover that $M + \mathbf{v}\mathbf{v}^\T$ is a rank one update of the matrix $M$. Thus, it is known that $M + \mathbf{v}\mathbf{v}^\T$ will be invertible if $1 + \mathbf{v}^\T M^{-1} \mathbf{v} \ne 0$; see Lemma 1.1 of \citet{ding2007eigenvalues}. In that case, we know that the unique solution to this system is $\phi = (M + \mathbf{v}\mathbf{v}^\T)^{-1}\mathbf{v}$.

As we noted in the proof of Theorem~\ref{orderStatRand}, the steady-state distribution of an $M/G/\infty$ queue is $\mathrm{Pois}(\lambda \E{S})$ when the arrival rate is $\lambda$ and service distribution is equivalent to the random variables $S$. We can now note further that $\lambda \E{S}$ is the steady-state mean of such a queueing system. The distribution of $Q_i$ is then given by $\mathrm{Pois}(\E{Q_i})$ for each $i \in \{1, \dots, k\}$, and thus is equivalent across all sub-queues.
\end{proof}
\end{proposition}

For added motivation, we now consider the two dimensional case in the following remark.

\begin{remark}
If $k = 2$, $M$ and $\phi$ are scalars, given by
$$
M
=
\frac{\E{S}}{\E{S_{2,2}} - \E{S_{1,2}}}
,
\quad
\phi = \frac{\E{S_{2,2}} - \E{S_{1,2}}}{\E{S} + \E{S_{2,2}} - \E{S_{1,2}}}
.
$$
In this case, we can note that if $\PP{B = 1} = \phi$, then in steady-state the distribution of the workload in the system from the easier jobs from all batches will be equivalent to that of the harder jobs. If $\PP{B = 1} > \phi$ the number of harder jobs will stochastically dominate the number of easier jobs, and vice-versa is $\PP{B=1} < \phi$.
\end{remark}

This result implies if we have the ability to choose the probability of batch sizes, we can construct each of the sub-systems which are organized by the order statitics to have the same queue length distribution.  Thus, providing equal work to all of the queues.

\section{Conclusion and Final Remarks} \label{concSec}

In this paper, we have found parallels between infinite server queues with batch arrivals, sums of scaled Poisson random variables, and Hermite distributions. Moreover, we also connect the stochastic objects to analytic quantities and functions of external interest, such as the harmonic numbers, the exponential integral function, the Euler-Mascheroni constant, and the polylogarithm function. In addition to being interesting in their own right, these connections have helped us to specify exact forms of valuable quantities related to this queueing system, including generating functions for the queue and for the limit of the queue scaled by the batch size. Thus, we have gained both insight into the queue itself and perspective on the model's place in operations research and applied mathematics more broadly.

For this reason, we believe continued work on these fronts is merited. For example, while we have some intuition for the harmonic Hermite distribution discussed in Subsection~\ref{markovSS}, we have less of an understanding of the limiting distribution of the scaled queue in that subsection and extended for random batch sizes in Subsection~\ref{limitRandSubsec}. Having more knowledge of what distribution might produce a moment generating function comprised of exponential integral function. Finding such a distribution could not only teach us about this queueing system, it would also likely be worth studying entirely on its own.  Additionally, providing further connections of this distribution back to the harmonic numbers and the associated Hermite distribution would also be of interest, such as in the connection of the limiting moment generating function to the expected value of a harmonic number evaluated at a Poisson random variable that we remarked in Subsection~\ref{markovSS}. One could also consider control problems for the routing of arrivals to sub-systems, like what we discuss for the case of random batch sizes in Subsection~\ref{orderStatRandSubsec}.

For future expansions of this work into other areas of queueing, we can group the main themes of potential further investigations in three categories.  First, the extension of our batch model beyond infinite server queues to multi-server queues, queues with abandonment, and networks of infinite server queues, a la \citet{mandelbaum2007service, massey2013gaussian, engblom2014approximations, gurvich2013excursion, pender2014gram, daw2017new}.  It would be interesting to explore our limit theorems in these cases to understand the impact of having a finite number of servers.  Second, it would also be interesting to explore the impact of the batch arrivals in the context of queues with delayed information as in \citet{pender2017queues, pender2017strong, pender2018analysis}.  It would be of interest to know whether or not the batch arrivals would influence the Hopf bifurcations or oscillations that occur in the delayed information queues. Additionally, one could explore findings of this work, like the steady-state distribution representation or the batch scaling, in contexts where there is dependence among the service durations within each batch of arrivals, such as those studied in \citet{pang2012infinite, falin1994m}. Finally, we are particularly interested in studying the impact of batch arrivals in the context of self-exciting arrival processes such as Hawkes processes like in the work of \citet{gao2016functional, koops2017infinite, daw2018queues}.   We intend to pursue the ideas described here as well as other related concepts in our future work.

\section*{Acknowledgements}
We acknowledge the generous support of the National Science Foundation (NSF) for Jamol Pender's Career Award CMMI \# 1751975 and Andrew Daw's NSF Graduate Research Fellowship under grant DGE-1650441.

\bibliographystyle{plainnat}
\bibliography{Batch}
%
%
%

\end{document}